\input amstex
\documentstyle{amsppt} 
\magnification=\magstep1
\NoRunningHeads
\NoBlackBoxes
\topmatter
\title Orbit equivalence of topological Markov shifts and Cuntz-Krieger algebras\endtitle
\author Kengo Matsumoto
\endauthor
\affil Department of Mathematical Sciences \\
Yokohama City University \\
Seto 22-2, Kanazawa-ku, Yokohama 236-0027, JAPAN
\endaffil
\abstract{We will prove that one-sided topological Markov shifts
$(X_A,\sigma_A)$ and $(X_B,\sigma_B)$ for matrices $A$ and $B$ with 
entries in $\{0,1\}$ 
are topologically orbit equivalent 
if and only 
if there exists an isomorphism between 
the Cuntz-Krieger algebras ${\Cal O}_A$ and ${\Cal O}_B$ 
keeping their commutative 
$C^*$-subalgerbas $C(X_A)$ and $C(X_B)$.
It is also equivalent to the condition that there exists a homeomorphism
from
$X_A$ to $X_B$ intertwining their topological full groups.  
We will also study structure of the automorphisms of ${\Cal O}_A$ keeping the commutative 
$C^*$-algebra $C(X_A)$.
 }
\endabstract
\endtopmatter


\def\Zp{{ {\Bbb Z}_+ }}

\def\OA{{{\Cal O}_A}}
\def\OB{{{\Cal O}_B}}
\def\FA{{{\Cal F}_A}}
\def\FB{{{\Cal F}_B}}
\def\DA{{{\frak D}_A}}
\def\DB{{{\frak D}_B}}
\def\HA{{{\frak H}_A}}
\def\HB{{{\frak H}_B}}

\def\Homeo{{{\operatorname{Homeo}}}}
\def\Out{{{\operatorname{Out}}}}
\def\Aut{{{\operatorname{Aut}}}}
\def\Inn{{{\operatorname{Inn}}}}
\def\Ker{{{\operatorname{Ker}}}}

\def\id{{{\operatorname{id}}}}
\def\supp{{{\operatorname{supp}}}}


\heading 1.Introduction
\endheading
Study of orbit equivalence of ergodic finite measure  preserving transformations was initiated by H. Dye [D], [D2], who proved that any such tranformations are orbit equivalent.
W. Krieger [Kr] has proved that two ergodic non-singular transformations are orbit equivalent if and only if the associated von Neumann crossed producs are isomorphic.
 In topological setting, Giordano-Putnam-Skau [GPS],[GPS2] (cf.[HPS]) have proved that two Cantor minimal systems are strong orbit equivalent 
 if and only if the associated $C^*$-crossed products are isomorphic.
In more general setting, 
J. Tomiyama [To]  (cf. [BT], [To2] ) has proved that    
two topological free homeomorphisms
$(X,\phi)$ and $(Y,\psi)$ on compact Hausdorff spaces
are continuously orbit equivalent if and only if 
there exists an isomorphism between the associated  $C^*$-crossed products 
keeping their commutative $C^*$-subalgebras $C(X)$ and $C(Y)$.
He also proved that it is equivalent to the condition that 
there exists a homeomorphism $h : X \rightarrow Y$ such that
$ h $ preserves their topological full groups.

In this paper we will study relationship between 
 orbit structure of one-sided topological Markov shifts and 
 algebraic structure of the associated Cuntz-Krieger algebras.
Let 
$(X_A,\sigma_A)$ 
be the right one-sided topological Markov shift
defined by an $N\times N$ square matrix $A$ with entries in $\{0,1\}$,
where $\sigma_A$ denotes the shift transformation.
 The one-sided topological Markov shifts 
 are no longer homeomorphism in general 
and the Cuntz-Krieger algebras 
 can not be written as a crossed product by $\Bbb Z$  in natural way.
 Hence Giordano-Putnam-Skau and Tomiyama's method 
 can not apply to study one-sided topological Markov shifts and Cuntz-Krieger algebras.  
However similar  type theorems to theirs 
will be proved  in our setting
by using a representation of $\OA$ on a Hilbert space having its complete orthonormal basis consisting of all points of the shift space $X_A$.
Let $\DA$ be the $C^*$-subalgebra consisting of all diagonal elements 
of the canonical AF-algebra $\FA$ inside of $\OA$. 
It is naturally isomorphic to the commutative $C^*$-algebra 
$C(X_A)$ of all ${\Bbb C}$-valued continuous functions on $X_A$.
Let $[\sigma_A]_c$ be the topological full group of 
$(X_A,\sigma_A)$ 
whose elements consist of homeomorphisms $\tau$ on $X_A$ 
such that $\tau (x) $ 
is contained in the orbit $orb_{\sigma_A}(x)$ of $x$
by $\sigma_A$ for all $x \in X_A$, 
and its orbit cocycles are continuous. 
We say that $(X_A, \sigma_A)$ and $(X_B, \sigma_B)$ 
are topologically orbit equivalent if there exists a homeomorphism $h:X_A \longrightarrow X_B$ 
such that
$h(orb_{\sigma_A}(x)) = orb_{\sigma_B}(h(x))$ for $x \in X_A$
and their orbit cocycles are continuous.
We will prove the following theorem:
\proclaim{Theorem 1.1(Theorem 5.6)}
Let $A, B$ be square matrices with entries in $\{0,1\}$.
Assume that the matrices $A, B$ and their transposed matrices $A^t, B^t$ 
 all satisfy condition (I) in [CK]. 
Then the  following are equivalent:
 \roster
 \item There exists an isomorphism $\Phi: \OA \rightarrow  \OB$
 such that $\Phi(\DA) = \DB$.
\item
$(X_A, \sigma_A)$ and $(X_B, \sigma_B)$ are topologically orbit equivalent.
\item There exists a homeomorphism $h: X_A \rightarrow X_B$ such that
$h \circ [\sigma_A ]_c \circ h^{-1} = [\sigma_B ]_c$. 
\endroster
\endproclaim
To prove the above theorem, 
we study the normalizer $N(\OA,\DA)$ of $\DA$ in $\OA$,
that is defined as  the group of all unitaries 
$u \in \DA$ such that $u\DA u^* =\DA$. 
We denote by ${\Cal U}(\DA)$ the group of all unitaries
in $\DA$.
We will prove 
\proclaim{Theorem 1.2(Theorem 4.8)}
Let $A$ be a square matrix with entries in $\{0,1\}$
satisfying condition (I) in [CK].
Then there exists a short exact  sequnce:
$$
1 
\longrightarrow {\Cal U}(\DA)
\overset{\id}\to{\longrightarrow} N(\OA,\DA)
\overset{\tau}\to{\longrightarrow} [\sigma_A ]_c
\longrightarrow 1
$$
that splits.
\endproclaim
Let $\Aut(\OA,\DA)$
be the group of all automorphisms $\alpha$ of 
$\OA$ such that $\alpha(\DA) = \DA$.  
Denote by $\Inn(\OA,\DA)$ the subgroup of $\Aut(\OA,\DA)$ 
of inner automorphisms
on $\OA$.
We set 
$
\Out(\OA,\DA) 
$
the quotient group
$\Aut(\OA,\DA)/\Inn(\OA,\DA).$
We will further prove 
\proclaim{Theorem 1.3(Theorem 6.5)}
Assume that both  $A$ and  $A^t$  satisfy condition (I) in [CK]. 
Then there exist short exact sequeneces:
\roster
\item
$
1 
\longrightarrow Z_{\sigma_A}^1({\Cal U}(\DA))
\overset{\lambda}\to{\longrightarrow} \Aut(\OA,\DA)
\overset{\phi}\to{\longrightarrow} N([\sigma_A ]_c)
\longrightarrow 1,
$
\item
$
1 
\longrightarrow B_{\sigma_A}^1({\Cal U}(\DA))
\overset{\lambda}\to{\longrightarrow} \Inn(\OA,\DA)
\overset{\phi}\to{\longrightarrow} [\sigma_A ]_c
\longrightarrow 1,
$
\item
$
1 
\longrightarrow H_{\sigma_A}^1({\Cal U}(\DA))
\overset{\lambda}\to{\longrightarrow} \Out(\OA,\DA)
\overset{\phi}\to{\longrightarrow} N([\sigma_A ]_c) / [\sigma_A ]_c
\longrightarrow 1.
$
\endroster
They all split. 
Hence  
$
\Out(\OA,\DA)
$
is a  semi-direct product
$$
\Out(\OA,\DA)
= N([\sigma_A ]_c) / [\sigma_A ]_c \cdot H_{\sigma_A}^1({\Cal U}(\DA)).
$$ 
\endproclaim
where $N([\sigma_A]_c)$ denotes the normalizer subgroup of
$[\sigma_A]_c$ in the group $\Homeo(X_A)$ of all homeomorphisms on $X_A$,
and
$Z_{\sigma_A}^1({\Cal U}(\DA))$,
$B_{\sigma_A}^1({\Cal U}(\DA))$
and
$H_{\sigma_A}^1({\Cal U}(\DA))$
are the group of unitary one-cocycles for $\sigma_A$,
the subgroup of $Z_{\sigma_A}^1({\Cal U}(\DA))$
of coboundaries and the cohomology group
$Z_{\sigma_A}^1({\Cal U}(\DA))/ B_{\sigma_A}^1({\Cal U}(\DA))$
respectively.

Similar type theorems hold for the pair of the canonical AF-algebra $\FA$  
inside $\OA$ and its diagonal algebra $\DA$, that are studied in Section 7.

The results of this paper will be generalized to more general  
subshifts and the $C^*$-algebras associated with the subshifts considered in
[Ma1] (cf.[CM]) and [Ma3] in a forthcoming paper [Ma4]. 
Throughout the paper,
we denote by
$\Zp$ and $\Bbb N$ the set of nonnegative integers and the set of positive integers
respectively.
\medskip

The author would like to thank to Takeshi Katsura for his excellent lectures 
on topological dynamics and $C^*$-algebras in COE seminor May 2007 at 
the University of Tokyo. 

\heading 2. Preliminaries
\endheading
Let $A=[A(i,j)]_{i,j=1}^N$ 
be an $N\times N$ matrix with entries in $\{0,1\}$,
where $1< N \in {\Bbb N}$.
Throughout the paper, we always assume that 
$A$ satisfies condition (I) in the sense of Cuntz-Krieger [CK].
We denote by 
$X_A$ the shift space 
$$
X_A = \{ (x_n )_{n \in \Bbb N} \in \{1,\dots,N \}^{\Bbb N}
\mid
A(x_n,x_{n+1}) =1 \text{ for all } n \in {\Bbb N}
\}
$$
over $\{1,\dots,N\}$
of the right one-sided topological Markov shift for $A$.
It is a compact Hausdorff space in natural  product topology.
The shift transformation $\sigma_A$ on $X_A$ is defined by 
$\sigma_{A}((x_n)_{n \in \Bbb N})=(x_{n+1} )_{n \in \Bbb N}$.
It is a continuous surjective map on $X_A$.
The topological dynamical system 
$(X_A, \sigma_A)$ is called the topological Markov shift for $A$.
The condition (I) for $A$
is equivalent to the condition that $X_A$ is homeomorphic to a Cantor discontinuum.

A word $\mu = \mu_1 \cdots \mu_k$ for $\mu_i \in \{1,\dots,N\}$
is said to be admissible for $X_A$ 
if $\mu$ appears in somewhere in some element $x$ in $X_A$.
 We denote by 
$B_k(X_A)$ the set of all admissible words of length $k \in \Bbb N$.
For $k=0$ we denote by $B_0(X_A)$ the empty word $\emptyset$.
We set 
$B_*(X_A) = \cup_{k=0}^\infty B_k(X_A)$ 
the set of admissible words of $X_A$.
For $x = (x_n )_{n \in {\Bbb N}} \in X_A$ and positive integers $k,l$ with
$k\le l$, we put the word 
$x_{[k,l]} = (x_k,x_{k+1},\dots, x_l) \in B_{l-k +1}(X_A)$
and the right infinite sequence
$x_{[k,\infty)} =(x_k, x_{k+1}, \dots ) \in X_A$.

The Cuntz-Krieger algebra $\OA$ for the matrix $A$ has been defined by 
the universal $C^*$-algebra generated by 
$N$ partial isometries $S_1,\dots, S_N$ subject to the relations:
$$ 
\sum_{j=1}^N S_j S_j^* = 1, \qquad
S_i^* S_i = \sum_{j=1}^N A(i,j) S_jS_j^*, \quad i=1,\dots,N \quad([CK]).
$$ 
If $A$ satisfies condition (I), 
the algebra $\OA$ is the unique $C^*$-algebra 
subject to the above relations.
For a word $\mu=\mu_1\cdots \mu_k \in B_k(X_A)$, we denote by 
$S_\mu = S_{\mu_1} \cdots S_{\mu_k}$.
By the universality for the above relations, 
the correspondence
$S_i \longrightarrow e^{\sqrt{-1}t}S_i, i=1,\dots,N$ 
for
$
e^{\sqrt{-1}t} \in 
{\Bbb T}=\{ e^{\sqrt{-1}t} \mid t \in [0,2 \pi]\}
$
yields an action 
$
\rho: {\Bbb T} \rightarrow \Aut(\OA)
$
that is called the gauge action.
It is well-known that the fixed point algebra of $\OA$ under $\rho$
is the AF-algebra $\FA$ generated by elements 
$S_\mu S_\nu^*, \mu, \nu \in B(X_A)$ 
with $|\mu | = |\nu|$ ([CK]). 
Let ${\Cal F}_A^n$ be the $C^*$-subalgebra of $\FA$
generated by elements $S_\mu S_\nu^*, \mu, \nu \in B_n(X_A)$.
Hence 
${\Cal F}_A^{\text{alg}} = \cup_{n=1}^{\infty}{\Cal F}_A^n$
is a dense $*$-subalgebra of $\FA$.  
We denote by $E: \OA \rightarrow \FA$
the conditional expectation define by
$E(a) = \int_{\Bbb T}\rho_t(a) dt $ for $a\in \OA$.
Let $\DA$ be the $C^*$-subalgebra of $\FA$ consisting of all diagonal elements 
of $\FA$. 
It is generated by elements 
$S_\mu S_\mu^*, \mu \in B_*(X_A)$
and isomorphic to the commutative $C^*$-algebra 
$C(X_A)$ of all ${\Bbb C}$-valued continuous functions on 
$X_A$ through the correspondence
$ 
S_\mu S_\mu^* \in \DA \longleftrightarrow \chi_\mu \in C(X_A)
$
where
$\chi_\mu$ denotes the characteristic function 
for the cylinder set 
$
U_\mu = \{ (x_n )_{n \in \Bbb N} \in X_A 
\mid x_1 =\mu_1,\dots, x_k = \mu_k \}
$
for $\mu=\mu_1\cdots\mu_k \in B_k(X_A)$.  
We identify $C(X_A)$ with the subalgebra 
$\DA$ of $\OA$.
Then the following lemma is well-known and basic in our further discussions.
\proclaim{Lemma 2.1( [CK;Remark 2.18], cf.[Ma2;Lemma 3.1])}
The algebra $\DA$ is maximal abelian in $\OA$.
\endproclaim

In [To], [To2], Tomiyama has used structure of pure state extensions of 
 point evaluations of the underlying space to study
orbit structure of topological dynamical systems of homeomorphisms on compact Hausdorff spaces.
However for the Cuntz-Krieger algebras, 
it has not been clarified structure of pure state extensions of 
point evaluations of the underlying shift space.
Instead of point evaluations, 
we will use a representation of the Cuntz-Krieger algebra 
$\OA$ on a Hilbert space having the shift space $X_A$ 
as a complete orthonormal basis, as in the following way.   
Let ${\frak H}_A$ 
be the Hilbert space with its complete orthonormal system
$ e_x , x \in X_A$. 
The Hilbert space is not separable.
Consider the partial isometries
$T_i, i=1,\dots,N$ defined by
$$
T_i e_x = 
\cases
e_{i x} & \text{ if } i x \in X_A,\\
0 & \text{ otherwise, }  
\endcases
$$
where $ix$ is defined by
$ix = (i,x_1,x_2,\dots )$ for $x = (x_n)_{n \in \Bbb N} \in X_A$.
It is easy to see that 
they satisfy the conditions
$ 
\sum_{j=1}^N T_j T_j^* = 1, 
T_i^* T_i = \sum_{j=1}^N A(i,j) T_jT_j^*
$ for 
$
 i=1,\dots,N 
$
so that 
the correspondence $S_i \rightarrow T_i$ 
yields a faithful representation of $\OA$ on ${\frak H}_A$.
We regard the algebra $\OA$ as the $C^*$-algebra generated by 
$T_i,i=1,\dots,N$
 on the Hilbert space
${\frak H}_A$ by this representation,
and write $T_i$ as $S_i$ (cf. [Ma2;Lemma 4.1]). 
\heading 3. Topological Full groups of Markov shifts
\endheading
For $x = (x_n )_{n \in \Bbb N} \in X_A$,
the orbit $orb_{\sigma_A}(x)$ of $x$ is defined by
$$
orb_{\sigma_A}(x) 
= \cup_{k=0}^\infty \cup_{l=0}^\infty \sigma_A^{-k}(\sigma_A^l(x)) \subset X_A.
$$
Hence  
$ y =( y_n )_{n \in \Bbb N} \in X_A$ 
belongs to $orb_{\sigma_A}(x)$ 
if and only if
there exists an admissible word 
$\mu_1 \cdots \mu_k \in B_k(X_A)$ 
such that 
$$
y = (\mu_1,\dots, \mu_k, y_{l+1}, y_{l+2},\dots )
\qquad \text{for some } k, l \in \Zp.
$$
We denote by $\Homeo(X_A)$
the group of all homeomorphisms on $X_A$.
We define the full group 
$[\sigma_A]$ and the topological full group
$[\sigma_A]_c$
for $(X_A,\sigma_A)$
as in the following way.

\noindent
{\bf Definition.}
Let
$
[\sigma_A ] 
$ be the set of all homeomorphism
$\tau \in \Homeo(X_A)
$ 
such that
$
\tau(x) \in  orb_{\sigma_A}(x)
$ 
for all
$ 
x \in X_A.
$ 
We call $[\sigma_A]$ the ful group of $(X_A,\sigma_A)$.
Let $[\sigma_A ]_c$ be the set of all $\tau$ in $[\sigma_A]$ such that  
there exist continuous maps 
$k, l : X_A \rightarrow \Zp$ 
such that 
$$
\sigma_A^{k(x)}(\tau(x) )=\sigma_A^{l(x)}(x)
\quad\text{ 
for all } x \in X_A. \tag 3.1
$$
We call   
$[\sigma_A ]_c$ the topological full group for $(X_A,\sigma_A)$.
The maps $k,l$ above are called orbit cocycles for $\tau$,
and sometimes written as $k_\tau, l_\tau$ respectively.
We remark that the orbit cocyles are not necessarily uniquely determined  
for $\tau$.

\noindent
{\bf Examples.}
(i) Put 
$F = 
\bmatrix
1 & 1 \\
1 & 0\\
\endbmatrix.
$
Define 
$\tau \in \Homeo(X_F)$ by setting
$$
\tau(x_1,x_2,\dots, )
=
\cases
(2,1,1,x_4, x_5,\dots, ) & \text{ if } (x_1,x_2,x_3) = (1,1,1), \\
(1,1,1,x_4, x_5,\dots, ) & \text{ if } (x_1,x_2,x_3) = (2,1,1), \\
(x_1,x_2,x_3,x_4, x_5,\dots, ) & \text{ otherwise.} \\
\endcases
$$
Since 
$\sigma_F(\tau(x)) = \sigma_F(x)$ for all $ x \in X_F$,
by putting
$k(x) = l(x) = 1$ for all $ x \in X_F$,
one sees that 
$\tau$ belongs to
$[\sigma_F]_c$.

(ii) More generally, 
let $A$ be an $N \times N$ matrix  with entries in $\{0,1\}$.
For $i \in \{ 1,\dots,N\}$ and $p \in {\Bbb N}$, 
we put
$$
W_p(i) = \{ (\mu_1,\dots, \mu_p) \in B_p(X_A) \mid A(\mu_p,i) =1 \}.
$$
We denote by 
${\frak S}(W_p(i))$ the group of all permutations on the set
$W_p(i)$.
Put 
$
{\frak S}_p(A) 
= {\frak S}(W_p(1)) \times \cdots \times {\frak S}(W_p(N)).
$
Then for an $N$-family 
$s = (s_1,\dots,s_N) \in {\frak S}_p(A) $
of permutations defines a homeomorphism $\tau_s \in \Homeo(X_A)$ 
by setting
$$
\tau_s(x_1,\dots, x_p,x_{p+1}, \dots )
= (s_{x_{p+1}}(x_1,\dots, x_p),x_{p+1}, \dots ),
\qquad x \in X_A.
$$
It is easy to see that 
$
\tau_s(x) \in  orb_{\sigma_A}(x)
$ 
for all
$ 
x \in X_A
$ 
and satisfies (3.1) 
for
$k(x) = l(x) = p$ for all $ x \in X_A$.
Hence
$\tau_s $ yields an element of $[\sigma_A]_c$.

\medskip

Let $A$ be an arbitrary fixed $N \times N$ matrix  with entries in $\{0,1\}$
satisfying condition (I).
The following lemma is direct.
\proclaim{Lemma 3.1}
$[\sigma_A]$ is a subgroup of $\Homeo(X_A)$ 
and 
$[\sigma_A]_c$ is a subgroup of $[\sigma_A]$. 
\endproclaim
Although $\sigma_A$ does not belong to $[\sigma_A]$,
the following lemma shows that $\sigma_A$ locally belongs to
$[\sigma_A]_c$, and  
the group $[\sigma_A]_c$ 
is not trivial in any case.
\proclaim{Lemma 3.2}
Assume that both the matrices $A$ and $A^t$ satisfy condition (I).
For any $\mu \in B_2(X_A)$,
there exists $\tau_\mu \in [\sigma_A]_c$
and continuous maps
$k_{\tau_\mu}, l_{\tau_\mu}: X_A \rightarrow \Zp$
such that  
$$
\cases 
\sigma_A^{k_{\tau_\mu}(x)}(\tau_\mu(x))
  = \sigma_A^{l_{\tau_\mu}(x)}(x)\qquad & \text{ for } x \in X_A,\\
\tau_\mu(y ) =\sigma_A(y) \qquad & \text{ for } y \in U_\mu,\\ 
k_{\tau_\mu}(y)  =0, \quad  l_{\tau_\mu}(y) =1 
\qquad & \text{ for } y \in U_\mu.\\
\endcases
\tag 3.2
$$
\endproclaim
\demo{Proof}
For 
$\mu = (\mu_1,\mu_2) \in B_2(X_A)$,
we have two cases.

 Case 1: $\mu_1= \mu_2$.

Put $a = \mu_1 = \mu_2$.
Since $A^t$ satisfies condition (I),
there exists $b_1 \in \{1,\dots, N\}$ 
such that 
$b_1 \ne a$ and $A(b_1,a) =1$.
Put
$\{b_1,\dots, b_{N-1} \} 
= \{1,\dots, N \} \backslash \{a\}$.
Let
$\{ b_{i_1}, \dots, b_{i_M}\} $
be the set of all elements of 
$\{b_1,\dots,b_{N-1}\}$ 
satisfying 
$A(a,b_{i_1}) = \cdots = A(a,b_{i_M}) =1$.
The set 
$\{b_{i_1}, \dots, b_{i_M}\}$ 
is nonempty because $A$ satisfies condition (I).
Define a homeomorphism
$\tau_\mu:X_A \rightarrow X_A$
by setting
$$
\tau_\mu(x) =
\cases
\sigma_A (x) \in U_a         & \text{ if } x \in U_{aa},\\ 
b_1 a b_{i_1}x_{[3,\infty)}\in U_{b_1 a b_{i_1}}
 & \text{ if } x= a b_{i_1}x_{[3,\infty)} 
\in U_{a b_{i_1}},\\
\vdots & \vdots \\
b_1 a b_{i_M}x_{[3,\infty)} \in U_{b_1 a b_{i_M}}& \text{ if } x= a b_{i_M}x_{[3,\infty)} 
\in U_{a b_{i_M}},\\
b_1 a a x_{[3,\infty)} \in U_{b_1 a a}& \text{ if } x= b_1 a x_{[3,\infty)} 
\in U_{b_1 a },\\
x & \text{ otherwise.}
\endcases
$$
We set
$$
k_{\tau_{\mu}}(x) =
\cases
0    & \text{ if } x \in U_{aa},\\ 
1  & \text{ if } x \in  U_{a b_{i_1}},\\
\vdots & \vdots \\
1 & \text{ if } x \in U_{a b_{i_M}},\\
2 & \text{ if } x \in U_{b_1 a },\\
0 & \text{ otherwise,}
\endcases
\qquad
l_{\tau_{\mu}}(x) =
\cases
1    & \text{ if } x \in U_{aa},\\ 
0  & \text{ if } x \in  U_{a b_{i_1}},\\
\vdots & \vdots \\
0 & \text{ if } x \in U_{a b_{i_M}},\\
1 & \text{ if } x \in U_{b_1 a },\\
0 & \text{ otherwise}
\endcases
$$
so that
$$
\sigma_A^{k_{\tau_\mu}(x)}(\tau_\mu(x))
  = \sigma_A^{l_{\tau_\mu}(x)}(x)\qquad \text{ for } x \in X_A.
$$
Hence 
$\tau_\mu \in [\sigma_A]_c$
and 
$\tau_\mu(y ) =\sigma_A(y),
k_{\tau_\mu}(y)  =0, 
l_{\tau_\mu}(y) =1
$
for
$ y \in U_\mu = U_{aa}.$

Case 2: $\mu_1 \ne \mu_2$.

Put $a = \mu_1, b = \mu_2$.
Define a homeomorphism
$\tau_\mu:X_A \rightarrow X_A$
by setting
$$
\tau_\mu(x) =
\cases
\sigma_A (x) \in U_b         & \text{ if } x \in U_{ab},\\ 
b x \in U_{ab}               & \text{ if } x \in U_b, \\
x & \text{ otherwise.}
\endcases
$$
We set
$$
k_{\tau_{\mu}}(x) =
\cases
0    & \text{ if } x \in U_{ab},\\ 
1  & \text{ if } x \in  U_{b},\\
0 & \text{ otherwise,}
\endcases
\qquad
l_{\tau_{\mu}}(x) =
\cases
1    & \text{ if } x \in U_{ab},\\ 
0  & \text{ if } x \in  U_{b},\\
0 & \text{ otherwise}
\endcases
$$
so that
$$
\sigma_A^{k_{\tau_\mu}(x)}(\tau_\mu(x))
  = \sigma_A^{l_{\tau_\mu}(x)}(x)\qquad \text{ for } x \in X_A.
$$
Hence 
$\tau_\mu \in [\sigma_A]_c$
and 
$\tau_\mu(y ) =\sigma_A(y),
k_{\tau_\mu}(y)  =0, 
l_{\tau_\mu}(y) =1
$
for
$ y \in U_\mu = U_{ab}.$
\qed
\enddemo
We remark that this lemma holds for any word $\mu$
with any length $|\mu| \ge 2$
by a similar argument to the above.

\proclaim{Lemma 3.3}
For $x = (x_n)_{n \in \Bbb N}\in X_A$ and
$j\in \{1,\dots,N\}$ 
with 
$jx = (j,x_1,x_2,\dots ) \in X_A$,
there exists 
$\tau \in [\sigma_A]_c$ 
such that $\tau(x) = jx$.
\endproclaim
\demo{Proof}
If $x= j^{\infty} = (j,j,\dots)$,
we may choose $\id$ as $\tau$.
If $x\ne j^{\infty}$, 
there exists
$k \in \Bbb N$ and $i \in \{ 1,\dots,N\}$ with $i \ne j$
such that
$x_n = j$ for $1 \le n \le k-1$
and $x_k = i$. 
Put
$\mu = (\undersetbrace{k-1}\to{j,\dots, j},i) \in B_k(X_A)$
and
$\nu =  (\undersetbrace{k}\to{j,\dots, j},i) = j \mu \in B_{k+1}(X_A)$
so that $x \in U_\mu$.
Define
$\tau: X_A \rightarrow X_A$ by setting
$$
\tau(y_1,y_2,y_3,\dots)
=
\cases
(j,y_1,y_2, \dots ) & \text{ if } y \in U_\mu,\\
(y_2,y_3,y_4,  \dots ) & \text{ if } y \in U_\nu,\\
(y_1,y_2, y_3, \dots ) & \text{ otherwise}.
\endcases
$$
Since
$U_\mu \cap U_\nu = \emptyset$,
one knows that $\tau: X_A \rightarrow X_A$
yields an element of $[\sigma_A]_c$.  
\qed
\enddemo

\proclaim{Lemma 3.4}
 Put
$[\sigma_A]_c (x) = \{ \tau(x) \in X_A \mid \tau \in [\sigma_A]_c \}$
for $x \in X_A$.
Then we have
$$
[\sigma_A]_c (x) = orb_{\sigma_A}(x)\qquad
\text{ for } x \in X_A.
$$
\endproclaim
\demo{Proof}
For any $\tau \in [\sigma_A]$, 
there exist
continuous maps
$k, l :X_A \rightarrow \Zp$
such that 
$\tau(x) = (\mu_1(x), \dots,\mu_{k(x)}(x), x_{l(x) +1},x_{l(x) +2}, \dots )
$ for some
$(\mu_1(x), \dots,\mu_{k(x)}(x)) \in B_{k(x)}(X_A)$
so that 
$\tau(x) \in orb(x)$ is clear,
and hence
$
[\sigma_A]_c (x) \subset orb_{\sigma_A}(x).
$

For the other inclusion relations,
by the previous lemmas,
for $x = (x_n )_{n \in \Bbb N} \in X_A$
and $j =\{1,\dots,N\}$ with $jx \in X_A$,
there exist
$\tau_1, \tau_2 \in [\sigma_A]_c$
such that
$$
\tau_1(x) = (j,x_1, x_2,\dots, ), \qquad
\tau_2(x) = (x_2, x_3,\dots, )
$$
so that
$$
[\sigma_A]_c(x) \ni (j,x_1, x_2,\dots, ),  (x_2, x_3,\dots, ).
$$  
Since $[\sigma_A]_c$ is a group, 
one sees that
$$
[\sigma_A]_c(x) \ni (\mu_1,\dots,\mu_k, x_{l+1},x_{l+2}, \dots, )
$$
for all $k,l\in \Zp$ 
and 
$\mu_1,\dots,\mu_k \in B_k(X_A)$
with
$(\mu_1,\dots,\mu_k, x_{l+1},x_{l+2},\dots ) \in X_A$.
Hence
$
[\sigma_A]_c (x) \supset orb_{\sigma_A}(x).
$
\qed
\enddemo

\heading 4. Full groups and normalizers  
\endheading
In this section, we will study 
the topological full group $[\sigma_A]_c$ and the normalizer 
$N(\OA,\DA)$.
We denote by ${\Cal U}(\OA), {\Cal U}(\DA)$ the groups of unitaries of
$\OA$ and $\DA$ respectively.
The normalizer $N(\OA,\DA)$ of $\DA$ in $\OA$ is defined by
$$
N(\OA,\DA) = \{ v \in {\Cal U}(\OA) \mid v \DA v^* = \DA \}.
$$
We will identify the algebra $C(X_A)$ with the subalgebra 
$\DA$ of $\OA$.  
We will first show the following proposition. 
\proclaim{Proposition 4.1}
For $\tau \in [\sigma_A]_c$, 
there exists a unitary $u_\tau \in N(\OA,\DA)$
such that
$$
Ad(u_\tau)(f) = f \circ \tau^{-1}
\quad
\text{ for } f \in \DA,
$$
and the correspondence 
$\tau \in [\sigma_A]_c \rightarrow u_\tau \in N(\OA,\DA)$  
is a homomorphism of group. 
\endproclaim
\demo{Proof}
Let the $C^*$-algebra $\OA$ be represented on the Hilbert space
${\frak H}_A$ with complete orthonormal basis 
$\{ e_x\mid x \in X_A\}$.
Then the generating partial isometries $S_i, i=1,\dots,N$ 
act on $\HA$ 
by 
 $S_i e_x = e_{ix}$ if $ix\in X_A$, otherwise $0$.
Since
 $\tau: X_A \rightarrow X_A$ is a homeomorphism,
 the operator $u_\tau$ on $\HA$ defined by 
 $$
 u_\tau e_x = e_{\tau(x)} \qquad \text{ for } x \in X_A
$$
yields a unitary on $\HA$.
We will prove that $u_\tau $ belongs to $\OA$.
As $\tau \in [\sigma_A]_c$, there exist
continuous maps
$l,k:X_A \rightarrow \Zp$
satisfying (3.1).
Since both $k(X_A)$ and $l(X_A)$ are finite sets 
of $\Zp$, there exist
$$
\tilde{k} = \max\{k(x) \mid x \in X_A \} \quad 
\text{ and }\quad
\tilde{l} = \max\{l(x) \mid x \in X_A \}\quad
\text{ in }
\Zp.
$$
Take 
$\mu(x) = (\mu_1(x),\dots, \mu_{k(x)}(x)) \in B_{k(x)}(X_A)$
such that
$$
\tau(x) = (\mu_1(x),\dots, \mu_{k(x)}(x),x_{l(x) +1},x_{l(x) +2},x_{l(x) +3},\dots)
$$
Since the set of words
$$
\{ (\mu_1(x), \dots, \mu_{k(x)}(x)) \in B_{k(x)}(X_A) \mid x\in X_A\}
$$
is a finite subset of
$$
W_{\tilde{k}}(X_A) = \cup_{j=0}^{\tilde{k}}B_j(X_A),
$$
the map
$x \in X_A \longrightarrow (\mu_1(x), \dots, \mu_{k(x)}(x))\in W_{\tilde{k}}(X_A)
$
is continuous, 
where 
$W_{\tilde{k}}(X_A)$
is endowed with discrete topology.
For any word 
$\nu \in W_{\tilde{k}}(X_A)$ with
$\nu = \nu_1 \cdots \nu_j, j \le \tilde{k}$,
the set
$$
 E_\nu = \{x \in X_A \mid \mu_1(x) = \nu_1,\dots, \mu_{k(x)}(x)= \nu_j\}
$$
is clopen such taht 
$\cup_{\nu \in W_{\tilde{k}}(X_A)}E_\nu = X_A$.
For $0 \le n \le \tilde{l}$,
the set 
$$
F_n = \{ x \in X_A \mid l(x) = n \}
$$
is  clopen in $X_A$.
We set 
$$
\align
Q_\nu & = \chi_{E_\nu} \in \DA \quad \text{ for } \nu \in W_{\tilde{k}}(X_A),\\
P_n & = \chi_{F_n} \in \DA \quad \text{ for } 0 \le n \le \tilde{l}.
\endalign
$$
Since
$X_A$ is disjoint unions
$
X_A = \cup_{\nu \in W_{\tilde{k}}(X_A)}E_\nu 
= \cup_{n=0}^{\tilde{l}} F_n,
$
one has
$$
\sum_{\nu \in W_{\tilde{k}}(X_A)} Q_\nu = \sum_{n=0}^{\tilde{l}} P_n = 1.
$$
For $ x \in X_A$ and  
$\nu \in W_{\tilde{k}}$, $0 \le n \le \tilde{l}$,
one has 
$x \in E_\nu \cap F_n$ if and only if 
$e_{\tau(x)} = S_\nu e_{\sigma_A^n(x)}$,
so that
$$
e_{\tau(x)} = \sum_{n=0}^{\tilde{l}} \sum_{\nu \in W_{\tilde{k}}}
( S_\nu \sum_{\xi \in B_n(X_A)} S_\xi^*) P_n Q_\nu e_x
\qquad \text{ for } x \in X_A.
$$ 
Therefore we have
$$
u_\tau =  \sum_{n=0}^{\tilde{l}} \sum_{\nu \in W_{\tilde{k}}}
( S_\nu \sum_{\xi \in B_n(X_A)} S_\xi^*) P_n Q_\nu
$$
that belongs to the algebra $\OA$.
The equality
$$
Ad(u_\tau)(f) = f \circ \tau^{-1} \quad
\text{ for } f \in \DA.
$$   
is straightforward
\qed.
\enddemo
For $v \in N(\OA,\DA)$,
we put
$Ad(v)(a) = vav^*$ for $a \in \OA$.
Then
$Ad(v)$ induces an automorphism on both algebras $\OA$ and $\DA$.
Let $\tau_v$ denotes the homeomorphism on $X_A$
induced by $Ad(v): \DA \rightarrow \DA$ 
satisfying
$Ad(v)(f) = f \circ {\tau_v}^{-1}$
for $f \in \DA$.
We will prove that $\tau_v$ 
gives rise to an element of $[\sigma_A]_c$.
We fix $v\in N(\OA,\DA)$ for a while.
\proclaim{Lemma 4.2}
There exists a family
$v_m, m\in {\Bbb Z}$ of partial isometries in $\OA$
such that all but finitely many $v_m, m\in {\Bbb Z}$ are zero, 
and
\roster
\item
$ v = \sum_{m \in {\Bbb Z}} v_m :$ finite sum. 
\item $v_m^*v_m, v_mv_m^*$ are projections in $\DA$. 
\item $v_m^* v_{m'} = v_m v_{m'}^* = 0$ for $m \ne {m'}$.
\item $v_0 \in \FA$ and satisfies
$v_0 \DA v_0^* \subset \DA$ and $v_0^* \DA v_0 \subset \DA$.
\endroster
\endproclaim
\demo{Proof}
Put $g(t) = v^* \rho_t(v) \in \OA$ for $t \in {\Bbb T}$.
For $f \in \DA$, one has          
$$
\rho_t(v) f \rho_t(v)^* = \rho_t(v f v^* ) = v f v^*
$$
so that $v^* \rho_t(v)$ commutes with each element of $\DA$.
By Lemma 2.1, $g(t)$ belongs to the algebra $\DA$.
We put
$$
v_m = \int_{\Bbb T}\rho_t(v) e^{-\sqrt{-1}mt} dt, 
\qquad 
\hat{g}(m) = \int_{\Bbb T} g(t) e^{-\sqrt{-1}mt} dt
\quad \text{ for } m \in {\Bbb Z}.
$$
Then $ v_m = v \hat{g}(m)$.  
Since $ g(t) \in \DA$,
one has
$$
\align
g(t)^* = & \rho_t(v^* \rho_{-t}(v)) = g(-t),\\
g(t) g(s) =& v^* \rho_t(v) \rho_t(v^* \rho_s(v)) = g(t+s).
\endalign
$$ 
One then sees that  
$\hat{g}(m), m\in {\Bbb Z}$ 
are projections in $\DA$ 
such that 
$\hat{g}(m)\hat{g}({m'}) =0$ for $m\ne {m'}$. 
Regard $g(t)\in \DA$ 
as a function on $X_A$.
For $x \in X_A$,
one sees that
$| g(t)(x)|^2 = <g(t)e_x \mid g(t)e_x > =1$ 
so that 
 by the Parseval's identity 
$$
1 = \int_{\Bbb T} | g(t)(x)|^2 dt 
= \sum_{m \in {\Bbb Z}} | \int_{\Bbb T}g(t)(x) e^{-\sqrt{-1}mt}dt|^2
 = \sum_{m \in {\Bbb Z}} \| \hat{g}(m)(x) e_x \|^2. 
$$
Put $E_m = \supp(\hat{g}(m)), m \in {\Bbb Z}$.
By the above equality, 
one has 
$X_A = \cup_{m \in {\Bbb Z}}E_m$ 
and
$E_m \cap E_{m'} = \emptyset$ for $m \ne {m'}$.
By the compactness of $X_A$,
one knows that all but finitely many $E_m$ are empty,
and that their sum is $1$.
Then both elements  
$v_m^*v_m = \hat{g}(m)$ and $v_m v_m^* = v \hat{g}(m) v^*$
are projections in $\DA$.
Therefore the assertions (1), (2)  and (3) hold.
For the assertion (4), 
we have
$$
v_0 = v \hat{g}(0) = v\int_{\Bbb T} v^* \rho_t(v) dt = E(v) \in \FA
$$ 
so that
$$
v_0 \DA v_0^* = v \hat{g}(0) \DA \hat{g}(0)v^* \subset \DA,
\qquad
v_0^* \DA v_0 = \hat{g}(0) v^* \DA v \hat{g}(0) \subset \DA,
$$
because  $\hat{g}(0) \in \DA$.
Thus we complete the proof.
\qed
\enddemo
\proclaim{Lemma 4.3}
For a fixed $n \in {\Bbb N}$, 
there exist partial isometries
$ v_\mu , v_{-\mu} \in \FA$ for each
 $\mu \in B_n(X_A)$
satisfying the following conditions:
\roster
\item
$ v_n = \sum_{\mu \in B_n(X_A)} S_\mu v_{\mu}$
and
$ v_{-n} = \sum_{\mu \in B_n(X_A)}  v_{-\mu}S_\mu^*$.
\item
 $v_{\mu}^*v_\mu$, $S_\mu v_\mu v_\mu^*S_\mu^*$, 
 $ S_\mu v_{-\mu}^*v_{-\mu}S_\mu^* $ 
and $v_{-\mu} v_{-\mu}^*$ are projections in $\DA$ such that
$$
\align
v_n^* v_n & = \sum_{\mu \in B_n(X_A)} v_\mu^* v_\mu, \qquad
v_n v_n^* = \sum_{\mu \in B_n(X_A)} S_\mu v_\mu v_\mu^* S_\mu^*,\\
v_{-n}^* v_{-n} & = \sum_{\mu \in B_n(X_A)} S_\mu v_{-\mu}^* v_{-\mu}S_\mu^*, \qquad
v_{-n} v_{-n}^* = \sum_{\mu \in B_n(X_A)} v_{-\mu} v_{-\mu}^*.
\endalign
$$
\item
$v_{\mu}v_\nu^* =  v_{-\mu}^* v_{-\nu}=0$ 
for $ \mu, \nu \in B_n(X_A)$ with 
$\mu \ne \nu$.
\item
The algebras 
$v_\mu \DA v_\mu^*, v_\mu^* \DA v_\mu, v_{-\mu} \DA v_{-\mu}^*$ 
and $v_{-\mu}^* \DA v_{-\mu}$ are contained in $\DA$. 
\endroster
\endproclaim
\demo{Proof}
Put for $ \mu \in B_n(X_A)$,
$$
v_\mu = E(S_\mu^* v), \qquad v_{-\mu} =E(v S_\mu). 
$$
They belong to $\FA$ and satisfy
$S_\mu^* S_\mu v_\mu = v_\mu$ and 
$v_{-\mu}S_\mu^* S_\mu = v_{-\mu}$.
Then we have
$$
\align
S_\mu^*v_n 
=&  \int_{\Bbb T}S_\mu^*\rho_t(v)e^{-\sqrt{-1}nt}dt=E(S_\mu^*v) =v_\mu,\\
v_{-n}S_\mu 
=&  \int_{\Bbb T}\rho_t(v)e^{\sqrt{-1}nt}S_\mu dt =E(vS_\mu)=v_{-\mu}.
\endalign
$$
Hence we have
$ v_n = \sum_{\mu \in B_n(X_A)} S_\mu v_{\mu}$
and
 $ v_{-n} = \sum_{\mu \in B_n(X_A)} v_{-\mu}S_\mu^*$.
Thus (1) holds.
We then have
$$
\align
v_{\mu}^*v_\mu 
& = v_n^* S_\mu S_\mu^*v_n = \hat{g}(n)v^* S_\mu S_\mu^* v \hat{g}(n),\\
S_\mu v_\mu v_\mu^* S_\mu^*
& =  S_\mu S_\mu^*v_n v_n^* S_\mu S_\mu^* =S_\mu S_\mu^* v \hat{g}(n) v^* S_\mu S_\mu^*,\\ 
S_\mu v_{-\mu}^*v_{-\mu}S_\mu^*
& = S_\mu S_\mu^* v_{-n}^* v_{-n}S_\mu S_\mu^* = S_\mu S_\mu^* \hat{g}(-n) S_{\mu}S_\mu^*,\\
v_{-\mu} v_{-\mu}^*
& = v_{-n}S_\mu S_\mu^* v_{-n}^* = v \hat{g}(-n) S_\mu S_\mu^* \hat{g}(-n)v^*.
\endalign
$$
Since
$\hat{g}(n),\hat{g}(-n)$ are projections in $\DA$
and $v\DA v^* =\DA$,
the above elements are projections in $\DA$ so that (2) and (3)  hold.
Since 
$$
v_\mu = S_\mu^* v_n = S_\mu^* v \hat{g}(n), \qquad
v_{-\mu} = v_{-n} S_\mu = v \hat{g}(-n) S_\mu
$$
the assertion (4) is immediate.
\qed
\enddemo

Let $u \in \OA$ be a partial isometry satisfying 
$$
u \DA u^* \subset \DA, \qquad
u^* \DA u \subset \DA.
$$
Put the projections
$p_u = u^* u, \, q_u = u u^* \in \DA$
and clopen sets
$X_u = \supp(p_u), \, Y_u = \supp(q_u) \subset X_A$.
Then $Ad(u) : \DA p_u \rightarrow \DA q_u$ 
yields an isomorphism
and induces a homeomorphism
$h_u: X_u \rightarrow Y_u$ such that
$$
Ad(u)(g) = g \circ {h_u}^{-1} \in \DA q_u (=C(Y_u))\qquad \text{ for }
g \in \DA p_u(= C(X_u)).
$$
\proclaim{Lemma 4.4}
Keep the above notation.
For $x \in X_u$, put $ y = h_u(x) \in Y_u$.
Then we have
$$
\| S^*_{y_{[1,n]}} u S_{x_{[1,n]}} \| = 1
\quad \text{ for all } n \in {\Bbb N}. 
$$
\endproclaim
\demo{Prof}
Since
$$
\|  S^*_{y_{[1,n]}} u S_{x_{[1,n]}} \|^2
 = \| S^*_{y_{[1,n]}} u S_{x_{[1,n]}}S_{x_{[1,n]}}^* u^* S_{y_{[1,n]}} \| 
$$
and
$S^*_{y_{[1,n]}} u S_{x_{[1,n]}}S_{x_{[1,n]}}^* u^* S_{y_{[1,n]}}$
is a projection in $\DA$,
one sees that 
$
\| S^*_{y_{[1,n]}} u S_{x_{[1,n]}} \| = 1
$
for all 
$
n \in {\Bbb N}
$
or 
$
\| S^*_{y_{[1,n]}} u S_{x_{[1,n]}} \| = 0
$
for all $n > N_0$ for some $N_0$.
It suffices to show that
$S^*_{y_{[1,n]}} u S_{x_{[1,n]}} \ne 0$ for all
$n \in \Bbb N$.
One then sees that
$$
(S^*_{y_{[1,n]}} u S_{x_{[1,n]}}
S_{x_{[1,n]}}^* u^* S_{y_{[1,n]}} 
e_{{\sigma_A}^n(y)} \mid  e_{{\sigma_A}^n(y)})
= ( Ad(u)(\chi_{U_{x_{[1,n]}}})e_y \mid e_y), 
$$
where
$\chi_{U_{x_{[1,n]}}}$ denotes the characteristic function on $X_A$
for the cylinder set $U_{x_{[1,n]}}$ of the word $x_{[1,n]}$.
As 
$$
Ad(u)(\chi_{U_{x_{[1,n]}}})e_y  
= (\chi_{U_{x_{[1,n]}}}\circ h_u^{-1})(y) e_y
= \chi_{U_{x_{[1,n]}}}(x) e_y
=e_y,
$$
one obtains that
$$
(S^*_{y_{[1,n]}} u S_{x_{[1,n]}}
S_{x_{[1,n]}}^* u^* S_{y_{[1,n]}} 
e_{{\sigma_A}^n(y)} \mid  e_{{\sigma_A}^n(y)})
= (e_y \mid e_y ) =1
$$
so that
$
S^*_{y_{[1,n]}} u S_{x_{[1,n]}} \ne 0.
$
\qed
\enddemo
The following is key
\proclaim{Lemma 4.5} 
Keep the above situation.
Assume that $u \in \FA$.
Then  there exists $k \in \Bbb N$ such that
for  all 
$
x=(x_n )_{n \in {\Bbb N}} \in X_u
$
$$
y_n = x_n
\qquad
\text{ for all }n > k
$$
where $y = h_u(x)$.
\endproclaim
\demo{Proof}
Suppose that for any $k \in {\Bbb N}$
there exist $x \in X_u$ and $N >k$ such that
$y_N \ne x_N$.
Now $u \in \FA$ so that take $u_{m_0} \in {\Cal F}_A^{k_0}$
 for some $k_0$ such that
$ \| u - u_{m_0} \| < \frac{1}{2}.$
Take $ x \in X_u$ and $N_0 > k_0$ such as $y_{N_0} \ne x_{N_0}$.
Since $u_{m_0}$ 
belongs to ${\Cal F}_A^{N_0-1}$,
it is written 
as 
$$
u_{m_0} = \sum_{\xi, \eta \in B_{N_0-1}(X_A)} 
c_{\xi,\eta}S_\xi S_\eta^* \in {\Cal F}_A^{N_0-1}
\quad \text{ for some } c_{\xi, \eta} \in {\Bbb C}.
$$
Hence we have
$$
S_{y_{[1,N_0-1]}}^* u_{m_0}S_{x_{[1,N_0-1]}}
=c_{y_{[1,N_0-1]},x_{[1,N_0-1]}} S_{y_{[1,N_0-1]}}^* S_{y_{[1,N_0-1]}}
S_{x_{[1,N_0-1]}}^*S_{x_{[1,N_0-1]}}
$$
so that
$$
\align
 & S_{y_{[1,N_0]}}^* u_{m_0}S_{x_{[1,N_0]}} \\
=& c_{y_{[1,N_0-1]},x_{[1,N_0-1]}}
 S_{y_{N_0}}^* S_{y_{[1,N_0-1]}}^* S_{y_{[1,N_0-1]}}
S_{x_{[1,N_0-1]}}^*S_{x_{[1,N_0-1]}}
 S_{x_{N_0}}
= 0
\endalign
$$
because
 $y_{N_0} \ne x_{N_0}$.
Hence we have
$$
S_{y_{[1,n]}}^* u_{m_0}S_{x_{[1,n]}} =0
\qquad 
\text{ for } n > N_0.
$$
For $n > N_0$, 
it then follows that
$$
\| S_{y_{[1,n]}}^* uS_{x_{[1,n]}} \|
= \| S_{y_{[1,n]}}^* ( u - u_{m_0})S_{x_{[1,n]}} \| < \frac{1}{2}.
$$ 
This is a contradiction to the preceding lemma. 
\qed
\enddemo
Thus we have
\proclaim{Lemma 4.6} 
For a partial isometry $u \in \FA$
satisfying 
$$
u \DA u^* \subset \DA, \qquad u^* \DA u \subset \DA,
$$
there exists $k_u \in {\Bbb N}$
such that the homeomorphism
$h_u : \supp(u^* u) \rightarrow \supp(u u^* )$
defined by
$Ad(u)(g) = g \circ {h_u^{-1}}$ for $ g \in \DA u^*u$
satisfies the condition
$$
\sigma_A^{k_u}(h_u(x)) = \sigma_A^{k_u}(x) \qquad \text{ for } 
x \in \supp(u^* u).
$$
\endproclaim
Therefore  we have
\proclaim{Proposition 4.7}
For any $v \in N(\OA,\DA)$,
the homomorphism $\tau_v$ 
induced by the automorphism on $\DA$
defined by the restriction to $Ad(v)$ to $\DA$ 
gives rise to an element of the topological full group
$[\sigma_A]_c$.
\endproclaim
\demo{Proof}
For $v \in N(\OA,\DA)$,
let $v_m,\, m\in {\Bbb Z}$ 
be the partial isometries in $\OA$ as in Lemma 4.2.
Take $K \in {\Bbb N}$ such that
$v_m = 0$ for all $m \in {\Bbb Z} $ with $|m | > K$,
and hence
$v = \sum_{m=-K}^K v_m$.
One has
$$
Ad(v)(f) = 
\sum_{n=1}^K v_n f v_n^*  
+ v_0 f v_0^*  
+ \sum_{n=1}^K v_{-n} f v_{-n}^*
\qquad \text{ for }
f \in \DA.
$$
As $v_m^*v_m, v_m v_m^*$ are projections in $\DA$,
we may put clopen sets 
$$
X_A^{(m)} = \supp(v_m^*v_m),\quad
Y_A^{(m)} = \supp(v_m v_m^*) \quad
\text{ for } m\in {\Bbb Z} 
\text{ with } |m| \le K
$$   
in $X_A$ such that
$X_A$ is the disjoint unions
$
X_A = \cup_{|m| \le K} X_A^{(m)}
=\cup_{|m| \le K} Y_A^{(m)}.
$ 
Since $v_0 \in \FA$, by Lemma 4.6, there exists 
$k_0\in \Bbb N$
such that 
$$
\sigma_A^{k_0}(\tau_0(x)) = 
\sigma_A^{k_0}(x)
\qquad
\text{ for } x \in X_A^{(0)}, \tag 4.1
$$
where
$\tau_0: X_A^{(0)} \longrightarrow Y_A^{(0)}$
is the homeomorphism
satisfying
$Ad(v_0)(f) = f \circ \tau_0^{-1}$
for
$f \in \DA v_0^* v_0$.
For
$v_n, v_{-n}, 1\le n \le K$, 
by Lemma 4.3,
one has
$$
\align
v_n f v_n^* & 
= \sum_{\mu \in B_n(X_A)} S_\mu v_\mu f v_\mu^* S_\mu^*, \\
v_{-n} f v_{-n}^* & 
= \sum_{\mu \in B_n(X_A)}  v_{-\mu} S_\mu^* f S_\mu v_{-\mu}^*
\quad 
\text{ for } f \in \DA.
\endalign
$$
Put
$$
\align 
X_A^{(n,\mu)} = \supp(v_\mu^*v_\mu),& \qquad
X_A^{(-n,\mu)} = \supp(S_\mu v_{-\mu}^*v_{-\mu}S_\mu^*)\\
\intertext{and}
Y_A^{(n,\mu)} = \supp(S_\mu v_\mu v_\mu^* S_\mu^*), & \qquad
Y_A^{(-n,\mu)} = \supp( v_{-\mu} v_{-\mu}^*).
\endalign
$$
By Lemma 4.3, one sees that
$$
X_A^{(m)} = \cup_{\mu \in B_{|m|}(X_A)}X_A^{(m,\mu)} ,\qquad
Y_A^{(m)} = \cup_{\mu \in B_{|m|}(X_A)}Y_A^{(m,\mu)}
$$
for $| m | \le K$.
There exist homeomorphisms
$$
\tau_{(m,\mu)}: X_A^{(m,\mu)} \longrightarrow Y_A^{(m,\mu)}
\qquad
\text{ for }
m \in {\Bbb Z} 
\text{ with }
|m | \le K
$$
such that
$$
\align
Ad(S_\mu v_\mu)(f) 
& = f \circ \tau_{(n,\mu)}^{-1} \qquad \text{ for } f \in \DA v_\mu^* v_\mu,\\
Ad(v_{-\mu} S_\mu^*)(g) 
& = g \circ \tau_{(-n,\mu)}^{-1} \qquad \text{ for } g \in 
\DA S_\mu v_{-\mu}^* v_{-\mu}S_\mu^*
\endalign
$$
for $n \in {\Bbb N}$ with $1 \le n \le K$.
As
$v_\mu, v_{-\mu} \in \FA$,
there exist
$k_{(n,\mu)},  k_{(-n,\mu)}\in {\Bbb N}$
such that
$$
\align
\sigma_A^{k_{(n,\mu)}}(\tau_{(n,\mu)}(x))
= & \sigma_A^{k_{(n,\mu)}+n}(x)
\qquad
\text{ for } x \in X_A^{(n,\mu)},\\
\sigma_A^{k_{(n,\mu)}+n}(\tau_{(-n,\mu)}(x))
= & \sigma_A^{k_{(n,\mu)}}(x)
\qquad
\text{ for } x \in X_A^{(-n,\mu)}.
\endalign
$$
Since we have
$$
\tau_v(x) =
\cases
\tau_{(n,\mu)}(x ) & \quad \text{ for }x \in X_A^{(n,\mu)},\\
\tau_0(x ) & \quad \text{ for }x \in X_A^{(0)},\\
\tau_{(-n,\mu)}(x ) & \quad \text{ for }x \in X_A^{(-n,\mu)}
\endcases
$$
and $X_A$ is the disjoint unions 
$$
\align
X_A &  = X_A^{(0)}  
 \bigcup \cup_{1\le |m|\le K}\cup_{\mu \in B_{|m|}(X_A)} X_A^{(m,\mu)} \\
    & =  Y_A^{(0)}  
 \bigcup \cup_{1\le |m|\le K}\cup_{\mu \in B_{|m|}(X_A)} Y_A^{(m,\mu)} 
\endalign
$$
of clopen sets
$
X_A^{(0)},X_A^{(m,\mu)}
$
and 
$
Y_A^{(0)},Y_A^{(m,\mu)}
$
for
$
1\le |m|\le K, \mu \in B_{|m|}(X_A),
$
  we conclude 
$\tau_v \in [\sigma_A]_c$.
\qed
\enddemo
The unitaries ${\Cal U}(\DA)$ are naturally embedded into 
$N(\OA,\DA)$. We denote the embedding by $\id$.
For $v \in  N(\OA,\DA)$, 
the induced homemorphism on $X_A$ is denoted by $\tau_v$,
that gives rise to an element of $[\sigma_A]_c$
by the above proposition.
We then have
\proclaim{Theorem 4.8}
The sequence
$$
1 
\longrightarrow {\Cal U}(\DA)
\overset{\id}\to{\longrightarrow} N(\OA,\DA)
\overset{\tau}\to{\longrightarrow} [\sigma_A ]_c
\longrightarrow 1
$$
is exact and  splits.
\endproclaim
\demo{Proof}
By Proposition 4.7, the map
$\tau: v \in N(\OA,\DA)
\longrightarrow 
\tau_v \in [\sigma_A ]_c
$
defines a homomorphism.
It is   surjective
by Proposition 4.1.
Suppose that $\tau_v =\id $ on $\DA$ 
for some $v \in N(\OA,\DA)$.
This means that $Ad(v) =\id $
on
$\DA$.
Hence $v$ commutes 
with all of elements of $\DA$.
By Lemma 2.1, $v$ belongs to $\DA$.
Therefore the sequence is exact.
As in Proposition 4.1,
for $\tau \in [\sigma_A]_c$,
the unitary $u_\tau$ defined by setting
$u_\tau e_x = e_{\tau(x)}, \, x \in X_A$ 
gives rise to a section of the exact sequence.
Hence the sequence splits.
\qed
\enddemo
\heading 5. Orbit equivalence
\endheading

\noindent
{\bf Definition.}
Let $(X_A, \sigma_A)$ and $(X_B,\sigma_B)$ be topological Markov shifts.
If there exists a homeomorphism 
$h:X_A \rightarrow X_B$ such that 
$h(orb_{\sigma_A}(x)) = orb_{\sigma_B}(h(x))$ for $x \in X_A$,
then 
$(X_A, \sigma_A)$ and $(X_B,\sigma_B)$ are said to be orbit equivalent.
In this case, for $x \in X_A$,
$h(\sigma_A(x)) \in orb_{\sigma_B}(h(x))$
so that 
$h(\sigma_A(x)) \in 
\cup_{k=0}^\infty \cup_{l=0}^\infty \sigma_B^{-k}\sigma_B^{l}(h(x))$.
Hence
there exist $k_1,\, l_1:X_A \rightarrow \Zp$
such that
$\sigma_B^{k_1(x)} (h(\sigma_A(x))) = \sigma_B^{l_1(x)}(h(x)) $.
Similarly 
there exist $k_2,\, l_2:X_B \rightarrow \Zp$
such that
$\sigma_A^{k_2(y)} (h^{-1}(\sigma_B(y))) = \sigma_A^{l_2(y)}(h^{-1}(y)) $.
 
We say that 
$(X_A, \sigma_A)$ and $(X_B,\sigma_B)$  
are
 {\it topologically orbit equivalent}\,\
if there exists a homeomorphism 
$h:X_A \rightarrow X_B$ such that 
$h(orb_{\sigma_A}(x)) = orb_{\sigma_B}(h(x))$ for $x \in X_A$
and
there exist 
continuous maps
$k_1,\, l_1:X_A \longrightarrow \Zp$
and
$k_2,\, l_2:X_B \longrightarrow \Zp$
such that
$$
\sigma_B^{k_1(x)} (h(\sigma_A(x))) = \sigma_B^{l_1(x)}(h(x)),
\qquad
\sigma_A^{k_2(y)} (h^{-1}(\sigma_B(y))) = \sigma_A^{l_2(y)}(h^{-1}(y)) 
\tag 5.1
$$
for $x \in X_A$ and $y \in X_B$.

\noindent
{\bf Example.}
Let 
$A_{[2]} =
\bmatrix
1 & 1 \\
1 & 1 
\endbmatrix,
\,
F =
\bmatrix
1 & 1 \\
1 & 0 
\endbmatrix.
$
The subshift $X_F$ is the set of all
sequences $(x_n)_{n \in \Bbb N}$ of $1,2 $ such that the word $(2,2)$ 
is forbidden.   
Define a homeomorphism
$h: X_F \longrightarrow X_{A_{[2]}}$
by substituing the word $(2,1)$ to the word $2$
from the leftmost in order such as
$$
\align
 h&(1,2,1,1,2,1,2,1,1,1,1,2,1,1,1,2,1,1,1,1,1,2,1,1,1,\dots) \\
= &(1,2,1,2,2,1,1,1,2,1,1,2,1,1,1,1,2,1,1,\dots) \in X_{A_{[2]}}.
\endalign
$$
Put for $i=1,2$
$$
\align 
U_{F,i}  
& =  
\{ x = (x_n)_{n \in \Bbb N}  \in X_{A_{[2]}} \mid x_1 =i\},\\
U_{A_{[2]},i} 
&=  \{ y =(y_n)_{n \in \Bbb N} \in X_{A_{[2]}} \mid y_1 =i\}.
\endalign
$$
By setting
$$
\cases
k_1(x) =0, \, l_1(x) =1 & \text{ for } x \in  U_{F,1}\\
k_1(x) =1, \, l_1(x) =1 & \text{ for } x \in  U_{F,2},
\endcases
\qquad
\cases
k_2(y) =0, \, l_2(y) =1 & \text{ for } y \in  U_{A_{[2]},1}\\
k_2(y) =0, \, l_2(y) =2 & \text{ for } y \in  U_{A_{[2]},2},
\endcases
$$
one knows that
$(X_F,\sigma_F)$
and 
$(X_{A_{[2]}}, \sigma_{A_{[2]}})$ 
are topologically orbit equivalent.

The following lemma is straightforward.
\proclaim{Lemma 5.1}
If $h: X_A \rightarrow X_B$ is a homeomorphism satisfying
$\sigma_B^{k(x)}(h(\sigma_A(x))) = \sigma_B^{l(x)}(h(x)), x \in X_A$
for some continous maps
$k,l:X_A \rightarrow \Zp$,
then by putting
$$
k^n(x) = \sum_{i=0}^{n-1}k(\sigma_A^i(x)),\qquad
l^n(x) = \sum_{i=0}^{n-1}l(\sigma_A^i(x)) 
$$
we have
$$
\sigma_B^{k^n(x)}(h(\sigma_A^n(x))) = \sigma_B^{l^n(x)}(h(x)),
\qquad x \in X_A, n \in \Bbb N.
$$
\endproclaim
\proclaim{Proposition 5.2}
If there exists a homeomorphism
$h:X_A \longrightarrow X_B$
such that
$ h \circ [\sigma_A]_c \circ h^{-1} = [\sigma_B]_c$,
then 
$(X_A, \sigma_A)$ and  
$(X_B, \sigma_B)$
are topologically orbit equivalent.
\endproclaim
\demo{Proof}
Assume that
$ h \circ [\sigma_A]_c \circ h^{-1} = [\sigma_B]_c$.
For any $y \in X_B$, put $x =h^{-1}(y)$
so that 
$h([\sigma_A]_c(x)) = [\sigma_B]_c(h(x)).$
By Lemma 3.3, we have
$[\sigma_A]_c(x) = orb_{\sigma_A}(x)$ and
$[\sigma_B]_c(h(x)) = orb_{\sigma_B}(h(x))$
so that
$h([\sigma_A]_c(x)) = orb_{\sigma_B}(h(x)).$

We will next show that there exist continuous cocycle functions
for $h$. 
By Lemma 3.2, 
For any $\mu \in B_2(X_A)$,
there exist $\tau_\mu \in [\sigma_A]_c$
and
$k_{\tau_\mu}, l_{\tau_\mu}: X_A \rightarrow \Zp$
satisfying (3.2).
Put 
$
\tau_h = h \circ \tau_\mu \circ h^{-1} \in h\circ [\sigma_A]_c \circ h^{-1} = [\sigma_B]_c.
$
For $ x \in U_\mu$, one has
$
h(\sigma_A(x)) = \tau_h(h(x)).
$
As 
$
\tau_h \in [\sigma_B]_c,
$
one may find 
$k_{\tau_h}^\mu, l_{\tau_h}^\mu: X_B \rightarrow \Zp$
such that
 $\sigma_B^{k_{\tau_h}^\mu(y)}(\tau_h(y)) = 
\sigma_B^{l_{\tau_h}^\mu(y)}(y).
$
For $y\in h(U_\mu)$, put
$x = h^{-1}(y)$
and hence
$$
\sigma_B^{k_{\tau_h}^\mu(h(x))}(h \circ \sigma_A(x)) = 
\sigma_B^{l_{\tau_h}^\mu(h(x))}(h(x))\qquad
\text{ for } x \in U_\mu.
$$
Let 
$\{ \mu^{(1)}, \dots, \mu^{(M)} \} $
be the set $B_2(X_A)$
of all admissible words of length $2$.
Define 
$k_1^h,l_1^h : X_A \rightarrow \Zp$
by setting
$$
k_1^h(x) = k_{\tau_h}^{\mu^{(i)}}(h(x)),\quad
l_1^h(x) = l_{\tau_h}^{\mu^{(i)}}(h(x))
\quad \text{ for } x \in U_{\mu^{(i)}}.
$$
They are continuous and satisfy 
$$
\sigma_B^{k_1^h(x)}(h \circ \sigma_A(x)) =
\sigma_B^{l_1^h(x)}(h(x)) 
\qquad
\text{ for } x \in X_A.
$$
Similarly there exist 
continious maps
$k_2^h,l_2^h : X_B \rightarrow \Zp$
such that
$$
\sigma_A^{k_2^h(y)}(h^{-1} \circ \sigma_B(y)) =
\sigma_A^{l_2^h(y)}(h^{-1}(y)) 
\qquad
\text{ for  } y \in X_B.
$$
Hence
$(X_A, \sigma_A)$ and  
$(X_B, \sigma_B)$
are topologically orbit equivalent.
\enddemo 
Conversely we have
\proclaim{Proposition 5.3}
If $(X_A, \sigma_A)$ and  
$(X_B, \sigma_B)$
are topologically orbit equivalent,
then there exists a homeomorphism
$h:X_A \longrightarrow X_B$
such that
$ h \circ [\sigma_A]_c \circ h^{-1} = [\sigma_B]_c$.
\endproclaim
\demo{Proof}
Suppose that 
there exists a homeomorphism 
$h:X_A \rightarrow X_B$ such that 
$h(orb_{\sigma_A}(x)) = orb_{\sigma_B}(h(x))$ for $x \in X_A$
and
there exist 
continuous maps
$k_1,\, l_1:X_A \rightarrow \Zp$
and
$k_2,\, l_2:X_B \rightarrow \Zp$
satisfying (5.1).
For $n \in {\Bbb N}$,
let 
$k_1^n, l_1^n:X_A \longrightarrow \Zp$
and 
$k_2^n, l_2^n:X_B \longrightarrow \Zp$
be continuous maps as in Lemma 5.1 such that
$$
\sigma_B^{k_1^n(x)} (h(\sigma_A^n(x)) = \sigma_B^{l_1^n(x)}(h(x)),
\qquad
\sigma_A^{k_2^n(y)} (h^{-1}(\sigma_B^n(y)) = \sigma_A^{l_2^n(y)}(h^{-1}(y)) 
\tag 5.2
$$
for $x \in X_A$ and $y \in X_B$.
For any $\tau \in [\sigma_A]_c$, 
there exist continuous maps:
$k_\tau, l_\tau : X_A \longrightarrow \Zp$
such that
$$
\sigma_A^{k_\tau(x)}(\tau(x)) = \sigma_A^{l_\tau(x)}(x)
\qquad
x \in X_A. \tag 5.3
$$
For $ y \in X_B$,
put
$x = h^{-1}(y)$.
We set
$m=k_\tau(x)\in \Bbb N$.
By (5.2) and (5.3), one has
$$
\sigma_B^{l_1^m (\tau(x))}(h(\tau(x))
=\sigma_B^{k_1^m(\tau(x))} (h(\sigma_A^m(\tau(x)))
=\sigma_B^{k_1^m(\tau(x))} (h(\sigma_A^{l_\tau(x)}(x))
$$
We set
$n=l_\tau(x)\in \Bbb N$.
By applying $\sigma_B^{k_1^n(x)}$ to the above equality,
one has by (5.2)
$$
\align
& \sigma_B^{k_1^n(x) +l_1^m (\tau(x))}(h(\tau(x))\\
=& \sigma_B^{k_1^m(\tau(x))}\sigma_B^{k_1^n(x)} (h(\sigma_A^n(x)))
=  \sigma_B^{k_1^m(\tau(x))}\sigma_B^{l_1^n(x)} (h(x))
=  \sigma_B^{k_1^m(\tau(x))+ l_1^n(x)} (h(x))
\endalign
$$
and hence
$$
\sigma_B^{k_1^n(x) +l_1^m (\tau(x))}(h\circ \tau \circ h^{-1}(y))
=\sigma_B^{k_1^m(\tau(x))+ l_1^n(x)} (h(x)).
$$
By putting
$$
\align
k_\tau^h(y) 
& = k_1^n(x) +l_1^m (\tau(x))
= k_1^{l_\tau( h^{-1}(y))}(h^{-1}(y)) 
+l_1^{k_\tau(h^{-1}(y)} (\tau(h^{-1}(y))),\\
l_\tau^h(y) 
& = k_1^m(\tau(x)) +l_1^n (x)
= k_1^{k_\tau( h^{-1}(y))}(\tau(h^{-1}(y))) 
+l_1^{l_\tau(h^{-1}(y)} (h^{-1}(y)),
\endalign
$$
one has
$$
\sigma_B^{k_\tau^h(y)}(h\circ \tau \circ h^{-1}(y))
=\sigma_B^{l_\tau^h(y)} (h(x))
\qquad 
\text{ for all }
y \in X_B
$$
so that
$
h\circ \tau \circ h^{-1} \in [\sigma_B]_c
$
and 
$
h \circ [\sigma_A]_c \circ h^{-1} \subset [\sigma_B]_c.
$
Similarly we have
$
h^{-1}\circ  [\sigma_B]_c \circ h \subset [\sigma_A]_c
$
and conclude 
$
h \circ [\sigma_A]_c \circ h^{-1} = [\sigma_B]_c.
$
\qed
\enddemo

\proclaim{Proposition 5.4}
If there exists an isomorphism $\Psi: \OA \longrightarrow \OB$
such that
$\Psi(\DA) = \DB$,
then there exists a homeomorphism
$h:X_A \longrightarrow X_B$
such that
$ h \circ [\sigma_A]_c \circ h^{-1} = [\sigma_B]_c$.
\endproclaim
\demo{Proof}
By Theorem 4.8,
there exists an isomorphism
$\widetilde{\Psi}: [\sigma_A]_c \longrightarrow [\sigma_B]_c$
of group
such that the following diagrams are commutative:
$$
\CD
1 @ >>>  {\Cal U}(\DA) @>\id>> N(\OA,\DA)
  @>\tau>> [\sigma_A ]_c @>>>1 \\
@. @VV{\Psi |_{{\Cal U}(\DA)}}V
@VV{\Psi}V 
@VV{\widetilde{\Psi}}V @. \\
1 @ >>>  U(\DB) @>\id>> N(\OB,\DB)
  @>\tau>> [\sigma_B ]_c @>>>1. \\
\endCD
$$
For any $v \in N(\OA,\DA)$,
put
$Ad(v)(f) = vfv^*$ 
for
$f \in \DA$.
We identify $\DA$ with $C(X_A)$ in natural way.
Let $\tau_v\in \Homeo(X_A)$ 
be the homeomorphism on $X_A$ 
satisfying
$Ad(v)(f) = f \circ \tau_v^{-1}$ for $f \in \DA$. 
The identity 
$$
\Psi \circ Ad(v) \circ \Psi^{-1} = Ad(\Psi(v)) \qquad 
\text{ for } v \in N(\OA,\DA) \tag 5.4
$$
holds. 
Let $h:X_A \longrightarrow X_B$ be the homeomorphism
satisfying
$\Psi(f) = f \circ h^{-1}$ for $f \in \DA$.
As 
$\Psi :
N(\OA,\DA) \longrightarrow N(\OB,\DB)$
is an isomorphism and
$\{ \tau_v \mid v \in N(\OA,\DA)\} = [\sigma_A]_c$,
(5.4) implies that 
$h \circ [\sigma_A]_c \circ h^{-1} = [\sigma_B]_c.$
\qed
\enddemo

\proclaim{Proposition 5.5}
If $(X_A, \sigma_A)$ and  
$(X_B, \sigma_B)$
are topologically orbit equivalent,
then 
there exists an isomorphism $\Psi:\OA \longrightarrow \OB$
such that
$\Psi(\DA) = \DB$.
\endproclaim
\demo{Proof}
The proof is essentially same as the proof of Proposition 4.1.
For the sake of completeness, 
we will give a complete proof for this proposition.
Let $h : X_A \rightarrow X_B$ 
be a homeomorphism giving rise to topologically orbit equivalence 
$(X_A, \sigma_A)$ and  
$(X_B, \sigma_B)$.
Represent 
$\OA$ on $\HA$
and 
$\OB$ on $\HB$
as usual respectively.
We will prove that
there exists a unitary
$u_h :{\frak H}_A \rightarrow {\frak H}_B$
such that
$$
Ad(u_h)(\OA) = \OB, \quad 
Ad(u_h)(\DA) = \DB
\quad \text{ and }
Ad(u_h)(f) = f \circ h^{-1},
\quad f \in \DA.
$$
Since $h : X_A \rightarrow X_B$ 
gives  rise to a topological orbit equivalence between
$(X_A, \sigma_A)$ and  
$(X_B, \sigma_B)$,
there exist
continuous maps
$k_1, l_1, k_2, l_2: X_A \rightarrow \Zp$
satisfying (5.1).
We denote by $e^A_x: x \in X_A$ and 
$e^B_x: x \in X_B$
 the complete orthonomal systems on 
$\HA$ and  $\HB$ coming from the shift spaces respectively.
Define the unitary
$u_h :\HA \longrightarrow \HB$
by setting
$u_h e_x^A = e_{h(x)}^B$ for $x \in X_A$.
We will first  prove that
$Ad(u_h)(\OA)= \OB$.
We denote by $S_i^A$ and $S_i^B$  
the canonical generating partial isometries for $S_i$ in $\OA$ and in $\OB$ 
respectively. 
For $ y \in X_B$, one has
$$
u_h S_i^A u_h^* e_y^B = 
\cases
e_{h(ih^{-1}(y))}^B & \text{ if }   ih^{-1}(y) \in X_A,\\
0 & \text{ otherwise. }
\endcases
$$
We set
$X_B^{(i)} = \{ y \in X_B \mid ih^{-1}(y) \in X_A \}$.
Put
$z = ih^{-1}(y) \in X_A$.
As $h(\sigma_A(z)) =y$,
by (5.1) one has
$h(z) \in \sigma_B^{-l_1(z)}(\sigma_B^{k_1(z)}(y))$.
Hence
$h(z) = (\mu_1,\dots,\mu_{l_1(z)}(z), y_{k_1(z)+1}, y_{k_2(z)+1}, \dots )$
for some
$\mu_1\cdots\mu_{l_1(z)}(z) \in B_{l_1(z)}(X_B)$.
Since
both the maps 
$k_1, l_1 : X_A \rightarrow \Zp$ and
the map
$y \longrightarrow  z = i h^{-1}(y)$ 
are continuous, there exist numbers
$$
\tilde{k}_1 = \max\{k_1(z) \mid y \in X_B^{(i)} \}, \qquad
\tilde{l}_1 = \max\{l_1(z) \mid y \in X_B^{(i)} \}.
$$ 
The set 
$
\{(\mu_1(z),\dots,\mu_{l_1(z)}(z)) \in B_{l_1(z)}(X_B) \mid y \in X_B^{(i)} \}
$
of words 
is a finite subset of
$W_{\tilde{l}_1}(X_B) = \cup_{j=0}^{\tilde{l}_1}B_j(X_B)$.
As
the map
$y \in X_B^{(i)} 
\rightarrow 
(\mu_1(z),\dots,\mu_{l_1(z)}(z)) \in W_{\tilde{l}_1}(X_B)
$
is continuous, where 
$W_{\tilde{l}_1}(X_B)$
is endowed with discrete topology,
for a word $\nu=\nu_1\cdots \nu_j \in W_{\tilde{l}_1}(X_B)$,
the set 
$$
E_\nu =\{ y \in X_B \mid \mu_1(z) = \nu_1,\dots,\mu_{l_1(z)}(z) =\nu_j \}
\quad \text{ where } z = ih^{-1}(y)
$$
is clopen in $X_B^{(i)}$.
For $0 \le n \le \tilde{k}_1$, 
we put
$$
F_n = \{ y \in X_B \mid k_1(z) = n \}
$$
that is clopen in $X_B$.
We set
$$
\align
Q_\nu & = \chi_{E_\nu} \in \DB \quad \text{ for } \nu \in 
W_{\tilde{l}_1}(X_B),\\  
P_n & = \chi_{F_n} \in \DB \quad \text{ for } n \in \Bbb N.
\endalign
$$
We put the projection 
$P^{(i)} = \chi_{X_B^{(i)}} \in \DB$.
Since we have disjoint unions
$$
X_B^{(i)} 
= \cup_{\nu \in W_{\tilde{l}_1}(X_B)}E_\nu 
= \cup_{n=0}^{\tilde{k}_1}F_n,
$$
we have
$$
P^{(i)} =\sum_{\nu \in W_{\tilde{l}_1}(X_B)}Q_\nu 
= \sum_{n=0}^{\tilde{k}_1} P_n.
$$
For $y \in X_B^{(i)}$ and
$ \nu \in W_{\tilde{l}_1}(X_B), 0 \le n \le \tilde{k}_1$,  
we have
$ y \in E_\nu \cap F_n$ if and only if
$h(ih^{-1}(y)) = \nu \sigma_B^n(y)$,
and the later condition is equivalent to the condition 
that
$e_{h(ih^{-1}(y))}^B = S_{\nu}^B e_{\sigma_B^n(y)}^B$.
Since
$ y \in E_\nu \cap F_n$ if and only if
$P_n Q_\nu e_y^B = e_y^B$ and
$
e_{\sigma_B^n(y)}^B = \sum_{\xi\in B_n(X_B)}{S_\xi^B}^* e_y^B,
$
we have
for $y \in X_B^{(i)}$,
$$
e_{h(ih^{-1}(y))}^B 
= \sum_{n=0}^{\tilde{k}_1} \sum_{\nu \in W_{\tilde{l}_1}(X_B)}
( S_{\nu}^B \sum_{\xi\in B_n(X_B)}{S_\xi^B}^*) P_n Q_\nu e_y^B.
$$
Hence 
$$
u_h S_i^A u_h^* e_y^B  
= \sum_{n=0}^{\tilde{k}_1} \sum_{\nu \in W_{\tilde{l}_1}(X_B)}
( S_{\nu}^B \sum_{\xi\in B_n(X_B)}{S_\xi^B}^*) P_n Q_\nu e_y^B
\qquad \text{ for } y \in X_B^{(i)}.
$$
Therefore
we have
$$
u_h S_i^A u_h^*   
= \sum_{n=0}^{\tilde{k}_1} \sum_{\nu \in W_{\tilde{l}_1}(X_B)}
( S_{\nu}^B \sum_{\xi\in B_n(X_B)}{S_\xi^B}^*) P_n Q_\nu P^{(i)}.
$$
As $P_n, Q_\nu, P^{(i)}$
are projections in $\DB$,
we have
$Ad(u_h) (S_i^A) \subset \OB$  
so that 
$Ad(u_h) (\OA) \subset \OB.$  
Since $u_h^* = u_{h^{-1}}$, we symmetrically have
$Ad(u_h^*) (\OB) \subset \OA$
so that 
$Ad(u_h) (\OA) = \OB.$

It is direct to see that 
$Ad(u_h) (f) = f \circ h^{-1}$ for $ f \in \DA$
so that we have
$Ad(u_h) (\DA) = \DB.$
\qed
\enddemo
Therefore we have
\proclaim{Theorem 5.6}
Let $A, B$ be square matrices with entries in $\{0,1\}$.
The following three assertions  are  equivalent:
 \roster
 \item There exists an isomorphism $\Phi: \OA \rightarrow  \OB$
 such that $\Phi(\DA) = \DB$.
\item
$(X_A, \sigma_A)$ and $(X_B, \sigma_B)$ are topologically orbit equivalent.
\item There exists a homeomorphism $h: X_A \rightarrow X_B$ such that
$h \circ [\sigma_A ]_c \circ h^{-1} = [\sigma_B ]_c$. 
\endroster
\endproclaim


\heading 6. Normalizers of the full groups and automorphisms of $\OA$
\endheading
In this section, we will study the normalizer subgroup
$$ 
N([\sigma_A]_c)
= \{ \varphi \in \Homeo(X_A) \mid 
\varphi \circ \tau \circ \varphi^{-1} \in
[\sigma_A]_c \text{ for all } \tau \in [\sigma_A]_c \}
$$  
of $[\sigma_A]_c$ in $\Homeo(X_A)$,
related to the automorphisms
$\Aut(\OA,\DA)$.
We set
$$
\align
N[\sigma_A] 
 = &\{ h\in \Homeo(X_A) 
\mid h(orb_{\sigma_A}(x))=orb_{\sigma_A}(h(x))\text{ for } x \in X_A\}, \\  
N_c[\sigma_A] 
 = &\{ h\in N[\sigma_A]) 
\mid 
\text{there exist continuous map }
k_1, l_1, k_2, l_2:X_A \rightarrow \Zp \\
& \text{ such that }
\sigma_A^{k_1(x)}(
h(\sigma_A(x)))=\sigma_A^{l_1(x)}(h(x)),\\
& \qquad \qquad \sigma_A^{k_2(x)}(
h^{-1}(\sigma_A(x)))=\sigma_A^{l_2(x)}(h^{-1}(x))
\text{ for } x \in X_A\} 
\endalign
$$
\proclaim{Lemma 6.1}
$N_c[\sigma_A]$ is a subgroup of $N[\sigma_A]$.
\endproclaim
\demo{Proof}
It suffices to show that
for 
$\varphi, \psi \in N_c[\sigma_A]$,
the composition $\psi \circ \varphi $ 
belongs to $N_c[\sigma_A]$. 
For $n \in \Bbb N$,
take continuous maps
$k_{1,\varphi}^n,\, 
l_{1,\varphi}^n,\,
k_{1,\psi}^n,\,
l_{1,\psi}^n: X_A \longrightarrow \Zp$
such that
$$
\align
\sigma_A^{k_{1,\varphi}^n(x)}(\varphi(\sigma_A^n(x)) &= 
\sigma_A^{l_{1,\varphi}^n(x)}(\varphi(x)) \tag 6.1 \\ 
\sigma_A^{k_{1,\psi}^n(x)}(\psi(\sigma_A^n(x)) &= 
\sigma_A^{l_{1,\psi}^n(x)}(\psi(x)) \tag 6.2
\endalign
$$ 
As in Lemma 5.1, we write
$
k_{1,\varphi}^n, \,
l_{1,\varphi}^n, \,
k_{1,\psi}^n, \,
l_{1,\psi}^n
$
as
$
k_{\varphi}^n, \,
l_{\varphi}^n, \,
k_{\psi}^n, \,
l_{\psi}^n
$
respectively.
By (6.2) for $\varphi(\sigma_A(x))$ as $x$, 
one has
$$
\sigma_A^{k_{\psi}^n(\varphi(\sigma_A(x)))}
(\psi(\sigma_A^n(\varphi(\sigma_A(x))))) 
= 
\sigma_A^{l_{\psi}^n(\varphi(\sigma_A(x)))}(\psi(\varphi(\sigma_A(x)))).
$$
Put $n = k_\varphi^1(x)$ and $m = l_\varphi^1(x)$.
By (6.1) for $n=1$,
one has
$
\sigma_A^n(\varphi(\sigma_A(x))) 
= 
\sigma_A^m(\varphi(x)) 
$ 
so that
$$
\sigma_A^{k_{\psi}^n(\varphi(\sigma_A(x)))}
(\psi(\sigma_A^m(\varphi(x)))) 
= 
\sigma_A^{l_{\psi}^n(\varphi(\sigma_A(x)))}(\psi(\varphi(\sigma_A(x))))
$$
and hence
$$
\sigma_A^{k_{\psi}^n(\varphi(\sigma_A(x)))}
\sigma_A^{k_{\psi}^m(\varphi(x))}
(\psi(\sigma_A^m(\varphi(x)))) 
= 
\sigma_A^{k_{\psi}^m(\varphi(x)) + l_{\psi}^n(\varphi(\sigma_A(x)))}
(\psi(\varphi(\sigma_A(x)))).
$$
By (6.2) we have
$$
\sigma_A^{k_{\psi}^n(\varphi(\sigma_A(x)))}
\sigma_A^{l_{\psi}^m(\varphi(x))}
(\psi(\varphi(x))) 
= 
\sigma_A^{k_{\psi}^m(\varphi(x)) + l_{\psi}^n(\varphi(\sigma_A(x)))}
(\psi(\varphi(\sigma_A(x)))).
$$
By putting
$$
k_{\psi \varphi}(x) 
= k_{\psi}^m(\varphi(x)) + l_{\psi}^n(\varphi(\sigma_A(x))),
\quad
l_{\psi \varphi}(x) 
= l_{\psi}^m(\varphi(x)) + k_{\psi}^n(\varphi(\sigma_A(x))),
$$
where
$n = k_\varphi^1(x), m = l_\varphi^1(x).$
The maps 
$k_{\psi \varphi}, l_{\psi \varphi}:
X_A \longrightarrow \Zp$ are continuous and satisfy
$$
\sigma_A^{k_{\psi \varphi}(x)}(\psi \varphi(\sigma_A(x))) 
= \sigma_A^{l_{\psi \varphi}(x)}(\psi \varphi(x)).
$$
Similarly,  
we may find 
continuous maps
$k_{\varphi^{-1}\psi^{-1}}, l_{\varphi^{-1}\psi^{-1}}:
X_A \longrightarrow \Zp$ that satisfy
$$
\sigma_A^{k_{\varphi^{-1}\psi^{-1}}(x)}(\varphi^{-1}\psi^{-1}(\sigma_A(x))) 
= \sigma_A^{l_{\varphi^{-1}\psi^{-1}}(x)}(\varphi^{-1}\psi^{-1}(x)).
$$
Therefore we know
$\psi\circ  \varphi \in N_c[\sigma_A]$.
\qed
\enddemo

\proclaim{Lemma 6.2}
$N_c[\sigma_A] = N([\sigma_A]_c)$.
\endproclaim
\demo{Proof}
For $\varphi \in N_c[\sigma_A]$ and
$\tau \in [\sigma_A]_c$, we will first prove that
$\varphi \circ \tau \circ \varphi^{-1} \in  [\sigma_A]_c$.
For $n\in {\Bbb N}$,
take continuous maps 
$
k_\varphi^n, \,
 l_\varphi^n, \,
k_{\varphi^{-1}}^n, \,
 l_{\varphi^{-1}}^n :X_A \longrightarrow \Zp$
satisfying
$$
\align
\sigma_A^{k_{\varphi}^n(x)}(\varphi(\sigma_A^n(x)) &= 
\sigma_A^{l_{\varphi}^n(x)}(\varphi(x)) \tag 6.3 \\ 
\sigma_A^{k_{\varphi^{-1}}^n(x)}(\varphi^{-1}(\sigma_A^n(x)) &= 
\sigma_A^{l_{\varphi^{-1}}^n(x)}(\varphi^{-1}(x)) \tag 6.4
\endalign
$$ 
for all $x \in X_A$.
For $\tau \in [\sigma_A]_c$,
let 
$k_\tau:X_A \longrightarrow \Zp$
be a continuous map satisfying (3.1).
By (6.3) one has
$$
\sigma_A^{k_{\varphi}^n(\tau \varphi^{-1}(x))}
(\varphi(\sigma_A^n(\tau \varphi^{-1}(x))))
= 
\sigma_A^{l_{\varphi}^n(\tau \varphi^{-1}(x))}(\varphi(\tau \varphi^{-1}(x))).
$$  
Put
$y = \varphi^{-1}(x),  n=  k_\tau(y), m= l_\tau(y)$.
By (3.1), one has
$\sigma_A^n(\tau(y)) = \sigma_A^m(y))$
so that
we have
$$
\sigma_A^{l_{\varphi}^n(\tau (y))}(\varphi(\tau (y)))
=
\sigma_A^{k_{\varphi}^n(\tau (y))}
(\varphi(\sigma_A^m(y))).
$$
By applying 
$\sigma_A^{k_\varphi^m(y)}$ to the above equality,
one has by (6.3)
$$
 \sigma_A^{k_\varphi^m(y) + l_{\varphi}^n(\tau (y))}(\varphi(\tau (y)) 
= \sigma_A^{k_{\varphi}^n(\tau (y))}
\sigma_A^{k_\varphi^m(y)}(\varphi(\sigma_A^m(y)))
= \sigma_A^{k_{\varphi}^n(\tau (y))}
\sigma_A^{l_\varphi^m(y)}(\varphi(y)).
$$
Put
$$
k_{\varphi\tau \varphi^{-1}}(x) 
= k_{\varphi}^m(y) + l_{\varphi}^n(\tau(y)),
\quad
l_{\varphi\tau \varphi^{-1}}(x) 
= k_{\varphi}^n(\tau(y)) + l_{\varphi}^m(y),
 $$
where
$y = \varphi^{-1}(x), n = k_\tau(y), m = l_\tau(y).$
The maps 
$k_{\varphi\tau \varphi^{-1}}, l_{\varphi\tau \varphi^{-1}}:
X_A \longrightarrow \Zp$ are continuous and satisfy
$$
\sigma_A^{k_{\varphi\tau \varphi^{-1}}(x) }(\varphi(\tau (\varphi^{-1}(x))
= \sigma_A^{l_{\varphi\tau \varphi^{-1}}(x)}(x).
$$
Hence
$\varphi \circ \tau \circ \varphi^{-1} \in [\sigma_A]_c$
so that
$\varphi \in N([\sigma_A]_c)$.

We will next prove the other inclusion relation 
$N_c[\sigma_A]) \supset N([\sigma_A]_c).$
For $\varphi \in N([\sigma_A]_c)$ 
one has
$\varphi \circ [\sigma_A]_c \circ \varphi^{-1}(y) = [\sigma_A]_c(y)$
for all $y \in X_A$.
Put $x = \varphi^{-1}(y)$.
By Lemma 3.4, one sees that
$$
\varphi(orb_{\sigma_A}(x)) = [\sigma_A]_c(\varphi(x)) = 
 orb_{\sigma_A}(\varphi(x)).
$$
Let
$\{ \mu^{(1)}, \dots, \mu^{(M)} \} $ be the set $B_2(X_A)$
of  all admissible words of length $2$.
For each word $\mu^{(i)}$, Lemma 3.2 shows that there exists
$\tau_i \in [\sigma_A]_c$
and continuous maps
there exist 
$
k^{(i)},l^{(i)} : X_A \longrightarrow \Zp
$ 
such that
$$
\tau_i(y) = \sigma_A(y) \quad \text{ for } y \in U_{\mu^{(i)}},
\quad
\sigma_A^{k^{(i)}(z)}(\tau_i(z)) = \sigma_A^{l^{(i)}(z)}(z)
\quad \text{ for } z \in X_A.
$$
Put
$
\tau_\varphi^{(i)} = \varphi \circ \tau_i \circ \varphi^{-1}$
that belongs to $[\sigma_A]_c.
$
We then have
$$
\varphi \circ \sigma_A(y) = \tau_\varphi^{(i)} (\varphi(y)) 
\quad \text{ for } y \in U_{\mu^{(i)}}.
$$
As $\tau_\varphi^{(i)} \in [\sigma_A]_c$
we may find continuous maps
$
k_{\tau_\varphi^{(i)}},
 l_{\tau_\varphi^{(i)}} :X_A \longrightarrow \Zp
$ 
such that
$$
\sigma_A^{k_{\tau_\varphi^{(i)}}(z)}(\tau_\varphi^{(i)}(z)) 
= \sigma_A^{l_{\tau_\varphi^{(i)}}(z)}(z)
\quad \text{ for } z \in X_A.
$$
Hence we have
$$
\sigma_A^{k_{\tau_\varphi^{(i)}}(y)}(\varphi \circ \sigma_A(y))
=\sigma_A^{l_{\tau_\varphi^{(i)}}(y)}(\varphi(y))
\quad \text{ for } y \in U_{\mu^{(i)}}.
$$
Define
$k_1^\varphi, l_1^\varphi : X_A \longrightarrow \Zp$ by setting
$$
k_1^\varphi(y) = k_{\tau_\varphi^{(i)}}(y),\quad
l_1^\varphi(y) = l_{\tau_\varphi^{(i)}}(y)
\quad
\text{ for } y \in U_{\mu^{(i)}}.
$$
Since
$ U_{\mu^{(i)}}$ 
is clopen and
$X_A $ is disjoint union $ \cup_{i=1}^M U_{\mu^{(i)}}$,
the maps
$k_1^\varphi, l_1^\varphi$
are both continuous and satisfy
$$
\sigma_A^{k_1^\varphi(y)}(\varphi \circ \sigma_A(y))
=\sigma_A^{l_1^\varphi(y)}(\varphi(y))
\quad \text{ for } y \in X_A.
$$
Similarly we may find continuous maps
$k_2^\varphi, l_2^\varphi : X_A \longrightarrow \Zp$ 
that satisfy
$$
\sigma_A^{k_2^\varphi(y)}(\varphi^{-1} \circ \sigma_A(x))
=\sigma_A^{l_2^\varphi(y)}(\varphi^{-1}(x))
\quad \text{ for } x \in X_A,
$$
so that 
$\varphi \in N_c[\sigma_A]$.
Therefore we have
$N_c [\sigma_A] \supset N([\sigma_A]_c)$
and hence 
$N_c [\sigma_A] = N([\sigma_A]_c)$.
\qed
\enddemo

\proclaim{Proposition 6.3}
For a homeomorphism $h \in N_c([\sigma_A])$
there exists an automorphism
$\alpha_h \in \Aut(\OA,\DA)$ such that 
$\alpha_h(f) = f \circ h^{-1}$
for
$f \in \DA$, and
the correspondence
$h \in N_c([\sigma_A]) 
\rightarrow \alpha_h \in \Aut(\OA,\DA)$
is a homomorphism.
\endproclaim
\demo{Proof}
Since a homomorphism
$h \in N_c([\sigma_A])$ gives rise to a topological 
orbit equivalence on $(X_A,\sigma_A)$,
the assertion follows from Proposition 5.5 and its proof.
\qed
\enddemo
Conversely
for any automorphism
$\alpha \in \Aut(\OA,\DA)$, we denote by
$\phi_\alpha$ the homeomorphism on $X_A$ induced by the restriction
of $\alpha$ to $\DA$ such that
$\alpha(f) = f \circ \phi_\alpha^{-1}$ for $f \in \DA$. 
We then have
\proclaim{Proposition 6.4}
$\phi_\alpha$ belongs to $N([\sigma_A]_c)$. 
\endproclaim
\demo{Proof}
For $\tau \in [\sigma_A]_c$, 
we will prove that
$\phi_\alpha\circ \tau \circ \phi_\alpha^{-1} \in [\sigma_A]_c.$
Let $u_\tau \in N(\OA,\DA)$ be the unitary constructed in Proposition 4.1,
such that 
$Ad(u_\tau)(f) = f \circ \tau^{-1}$ for $f \in \DA$.
Since
$
Ad(\alpha(u_\tau))
=\alpha \circ Ad(u_\tau) \circ \alpha^{-1}
$
on $\OA$,
the condition
$\alpha(\DA) = \DA$
implies
$\alpha(u_\tau)\in N(\OA,\DA)$.
One then sees that
$$
Ad(\alpha(u_\tau))(f) 
=\alpha \circ Ad(u_\tau) \circ \alpha^{-1}(f)
=f \circ (\phi_\alpha \circ \tau^{-1} \circ \phi_\alpha^{-1}).
$$
Since the homeomorphism
$\tau_{\alpha(u_\tau)}$ defined by $\alpha(u_\tau)\in N(\OA,\DA)$
belongs to
$[\sigma_A]_c$ 
and satisfies 
$
Ad(\alpha(u_\tau))(f) 
= f \circ \tau_{\alpha(u_\tau)}^{-1}$,
one concludes that
$$
\tau_{\alpha(u_\tau)}^{-1} 
= (\phi_\alpha \circ \tau^{-1} \circ \phi_\alpha^{-1})^{-1}
= \phi_\alpha \circ \tau \circ \phi_\alpha^{-1}
$$
that belongs to
$[\sigma_A]_c$.
\qed
\enddemo

We denote by $\varphi_A:\DA \rightarrow \DA$ the homomorphism
defined by 
$$
\varphi_A(a) = \sum_{i=1}^N S_i a S_i^* \qquad
\text{ for } a \in \DA.
$$
In regarding $\DA$ as $C(X_A)$ as usual,
one sees
$\varphi_A(f) = f\circ \sigma_A$ for $f \in C(X_A)$.
A unitary one-cocycle 
$U:\Zp \rightarrow {\Cal U}(\DA)$
for $\varphi_A$
is a ${\Cal U}(\DA)$-valued 
 function on $\Zp$ satisfying
$$
U(k+l) = U(k) \varphi_A^k(U(l), \qquad k,l\in \Zp \quad (cf. [Ma2]).
$$
Let
$Z_{\sigma_A}^1({\Cal U}(\DA))$
be the set of all unitary one-cocycles for $\varphi_A$,
that is an abelian group in natural way.
As in [Ma2] (cf. [C2], [KT]),
For $U \in Z_{\sigma_A}^1({\Cal U}(\DA))$,
put
$$
\lambda(U)(S_\mu) = U(k) S_\mu
\qquad \text{ for } \mu \in B_k(X_A).
$$
Then  
$
\lambda(U)
$
gives rise to an automorphism of $\OA$
such that $\lambda(U)|_{\DA} = \id$.
We note that the correspondence
$U \in Z_{\sigma_A}^1({\Cal U}(\DA)) \rightarrow U(1) \in {\Cal U}(\DA)$
yeilds an isomorphism of abelian group, and hence we may identify
$Z_{\sigma_A}^1({\Cal U}(\DA))$ with ${\Cal U}(\DA)$.
By [Ma2; Lemma 4.8],
$$
\lambda:Z_{\sigma_A}^1({\Cal U}(\DA)) \rightarrow \Aut(\OA, \DA)
$$
is an injective homomorphism of group.

Let $V: \Zp \rightarrow {\Cal U}(\DA)$ be a ${\Cal U}(\DA)$-valued function on 
$\Zp$ satisfing 
$$
U(k) = v \varphi_A^k(v^*), \qquad k \in \Zp
$$
for some unitary
$v \in {\Cal U}(\DA)$.
Then $V$ is called a coboundary for $\varphi_A$.
Since
$$
V(k) \varphi_A^k(V(l)) = v \varphi_A^k(v^*)\varphi_A^k(v \varphi_A^l(v^*))=V(k+l)
$$
 a coboundary $V$ for $\varphi_A$ is a unitary one-cocycle for 
$\varphi_A$.
Let
$B_{\sigma_A}^1({\Cal U}(\DA))$
be the set of all coboundaries for $\varphi_A$.
It is easy to see that
$B_{\sigma_A}^1({\Cal U}(\DA))$ is a subgroup of
$Z_{\sigma_A}^1({\Cal U}(\DA))$.
We remark that 
if $U \in Z_{\sigma_A}^1({\Cal U}(\DA))$
satisfies
$U(1) = v \varphi_A(v^*)$
for some
$v \in {\Cal U}(\DA)$,
then 
$U(k) = v\varphi_A^k(v^*)$ for $k\in \Bbb N$,
and hence $U \in B_{\sigma_A}^1({\Cal U}(\DA))$.

Define 
$H_{\sigma_A}^1({\Cal U}(\DA))$
by the quotient group
$Z_{\sigma_A}^1({\Cal U}(\DA))/B_{\sigma_A}^1({\Cal U}(\DA))$,
that is called the cohomology group for
$\varphi_A$.

\proclaim{Theorem 6.5}
There exist short exact sequences: 
\roster
\item
$
1 
\longrightarrow Z_{\sigma_A}^1({\Cal U}(\DA))
\overset{\lambda}\to{\longrightarrow} \Aut(\OA,\DA)
\overset{\phi}\to{\longrightarrow} N([\sigma_A ]_c)
\longrightarrow 1,
$
\item
$
1 
\longrightarrow B_{\sigma_A}^1({\Cal U}(\DA))
\overset{\lambda}\to{\longrightarrow} \Inn(\OA,\DA)
\overset{\phi}\to{\longrightarrow} [\sigma_A ]_c
\longrightarrow 1,
$
\item
$
1 
\longrightarrow H_{\sigma_A}^1({\Cal U}(\DA))
\overset{\lambda}\to{\longrightarrow} \Out(\OA,\DA)
\overset{\phi}\to{\longrightarrow} N([\sigma_A ]_c) / [\sigma_A ]_c
\longrightarrow 1.
$
\endroster
They all split. Hence 
$
\Out(\OA,\DA)
$
is a semi-direct product 
$$
\Out(\OA,\DA)
= N([\sigma_A ]_c) / [\sigma_A ]_c \cdot H_{\sigma_A}^1({\Cal U}(\DA)).
$$ 
\endproclaim
\demo{Proof}
(1)
As $N([\sigma_A]_c) =N_c[\sigma_A]$ by Lemma 6.1,
Proposition 6.3 and Proposition 6.4  imply that the homomorphism
$\phi: \Aut(\OA,\DA) \longrightarrow  N([\sigma_A ]_c)$
is defined and surjective.
By [Ma2;Lemma 4.8],  
$\lambda: Z_{\sigma_A}^1({\Cal U}(\DA)) \longrightarrow \Aut(\OA,\DA)$
is injective.
Let $\alpha \in \Aut(\OA,\DA)$ be such that
$\Phi(\alpha) =\id$ and hence 
$\alpha|_{\DA} =\id$.
By [Ma2;Corollary 4.7],
$\alpha|_{\DA} =\id$ if and only if 
$\alpha = \lambda(U) $
for some $U \in Z_{\sigma_A}^1({\Cal U}(\DA)).$
Hence we have
$\Ker(\phi) = Z_{\sigma_A}^1({\Cal U}(\DA))$.
By Proposition 6.3, 
for 
 $\varphi \in N_c[\sigma_A]$, 
there exists an automorphism
$\alpha_{\varphi} \in \Aut(\OA,\DA)$,
that is of the form
$\alpha_{\varphi} = Ad(u_{\varphi}),$
where $u_{\varphi}:{\frak H}_A \longrightarrow {\frak H}_B$
is a unitary defined in the proof of Proposition 5.5.
It is clear to see that 
$\phi \circ \alpha_{\varphi} = \varphi$.
Hence the sequence splits.

(2)
Theorem 4.8 implies that
the homomorphism
$
\phi: \Inn(\OA,\DA)\longrightarrow [\sigma_A ]_c
$
is defined and surjective.
For 
$\alpha \in \Inn(\OA,\DA)$, 
take $v \in {\Cal U}(\OA)$ such that 
$\alpha = Ad(v)$.
Hence $v$ belongs to $N(\OA,\DA)$ such that ,
suppose that 
$\Phi(Ad(v)) = \id$ in $[\sigma_A]_c$.
By (1), there exists a cocycle
$U \in Z_{\sigma_A}^1({\Cal U}(\DA))$ such taht
$Ad(v) = \lambda(U)$.
By {Ma2;lemma 5.14],
one sees that $ v \in {\Cal U}(\DA)$ and
$U(1) = v \varphi_A(v^*)$.
Hence $U$ belongs to
$B_{\sigma_A}^1({\Cal U}(\DA))$.
As the sequence (1) splits, the section in (1) 
yields a section in (2). 
Hence (2) splits.
 
(3) The exact sequnce follows from (1) and (2),
that splits.
\qed
\enddemo

\heading 7. Orbit equivalence and AF-algebras
\endheading
In this section,
we will show that the discussions in the previous sections can be applied to 
the pair $(\FA,\DA)$ of the AF-algebra $\FA$ and its diagonal algebra $\DA$,
instead of the pair $(\OA,\DA)$ that we have studied.
For $x = (x_n )_{n \in \Bbb N} \in X_A$,
the {\it uniform orbit}\, $orb_{\sigma_A}[x]$ of $x$ is defined by
$$
orb_{\sigma_A}[x] 
= \cup_{k=0}^\infty \sigma_A^{-k}(\sigma_A^k(x)) \subset X_A.
$$
Hence  
$ y =( y_n )_{n \in \Bbb N} \in X_A$ 
belongs to $orb_{\sigma_A}[x]$ 
if and only if
there exist $k \in \Zp$ and an admissible word 
$\mu_1 \cdots \mu_k \in B_k(X_A)$ 
such that 
$$
y = (\mu_1,\dots, \mu_k, y_{k+1}, y_{k+2},\dots ).
$$
Let
$
[[\sigma_A] ] 
$ be the set of all homeomorphisms
$\tau \in \Homeo(X_A)
$ 
such that
$
\tau(x) \in  orb_{\sigma_A}[x]
$ 
for all
$ 
x \in X_A.
$ 
Let $[\sigma_A ]_{AF}$
 be the set of all $\tau$ in $[[\sigma_A]]$ 
such that  
there exists a continuous map 
$k : X_A \rightarrow \Zp$ 
such that 
$$
\sigma_A^{k(x)}(\tau(x) )=\sigma_A^{k(x)}(x)
\quad\text{ 
for all } x \in X_A. \tag 7.1
$$
We call   
$[\sigma_A ]_{AF}$ the AF-full group for $(X_A,\sigma_A)$.
As $X_A$ is compact, 
for a homeomorphism $\tau \in \Homeo(X_A)$,
$\tau$ belongs to $ [\sigma_A ]_{AF}$ 
if and only if 
there exists a constant number $k \in \Zp$ such that
$
\sigma_A^k(\tau(x) )=\sigma_A^k(x)
$ 
for all
$ x \in X_A$.
We set for $x \in X_A$,
$
[\sigma_A]_{AF}(x) = \{ \tau(x) \mid \tau \in [\sigma_A]_{AF}\}.
$
It is immediate to see that
$
[\sigma_A]_{AF}(x) = orb_{\sigma_A}[x].
$
Let
$N(\FA,\DA)$
be the normalizer  of $\DA$ in $\FA$,
that is defined as  the group of all unitaries 
$u \in \FA$ such that $u\DA u^* =\DA$. 
We note that the algebra $\DA$ is also maximal abelian in $\FA$.
By a similar manner to the proof of Proposition 4.1,
we have
\proclaim{Lemma 7.1}
For any $\tau \in [\sigma_A]_{AF}$, 
there exists a unitary $u_\tau \in N(\FA,\DA)$
such that
$$
Ad(u_\tau)(f) = f \circ \tau^{-1}
\quad
\text{ for } f \in \DA,
$$
and the correspondence 
$\tau \in [\sigma_A]_{AF} \rightarrow u_\tau \in N(\FA,\DA)$  
is a homomorphism of group. 
\endproclaim
 By Lemma 4.5 we have
 \proclaim{Lemma 7.2}
 Let $u \in N(\FA,\DA)$ be a unitary. 
Let $h_u\in \Homeo(X_A)$ be the homeomorphism on $X_A$ induced 
by the restriction of $Ad(u)$ to $\DA$ such that
$Ad(u)(f) = f \circ h_u^{-1}$ for $f \in \DA$.
Then there exists a number $k \in {\Bbb N}$ 
such that
$
\sigma_A^k(h_u(x)) = \sigma_A^k(x)
$
for $x \in X_A$.  
Namely $h_u \in  [\sigma_A]_{AF}.$
 \endproclaim
Therefore by a similar proof of Theorem 4.8,
we have 
\proclaim{Prposition 7.3}
There exists a short exact  sequnce:
$$
1 
\longrightarrow {\Cal U}(\DA)
\overset{\id}\to{\longrightarrow} N(\FA,\DA)
\overset{\tau}\to{\longrightarrow} [\sigma_A ]_{AF}
\longrightarrow 1
$$
that splits.
\endproclaim

We say that 
$(X_A, \sigma_A)$ and $(X_B,\sigma_B)$  
are
 {\it uniformly orbit equivalent}\,\
if there exists a homeomorphism 
$h:X_A \rightarrow X_B$ such that 
$h(orb_{\sigma_A}[x]) = orb_{\sigma_B}[h(x)]$ for $x \in X_A$
and
there exist constant numbers
$k_1,\, k_2\in \Zp$
such that
$$
\sigma_B^{k_1} (h(\sigma_A(x)) = \sigma_B^{k_1}(h(x)),
\qquad
\sigma_A^{k_2} (h^{-1}(\sigma_B(y)) = \sigma_A^{k_2}(h^{-1}(y)) 
$$
for $x \in X_A$ and $y \in X_B$.
Then by a completely similar manner to
the proof of Proposition 5.2, Proposition 5.3, Proposition 5.4, Proposition 5.5
and Theorem 5.6,
we have
\proclaim{Theorem 7.4}
The following three assertions  are  equivalent:
 \roster
 \item There exists an isomorphism $\Phi: \FA \rightarrow  \FB$
 such that $\Phi(\DA) = \DB$.
\item
$(X_A, \sigma_A)$ and $(X_B, \sigma_B)$ are uniformly orbit equivalent.
\item There exists a homeomorphism $h: X_A \rightarrow X_B$ such that
$h \circ [\sigma_A ]_{AF} \circ h^{-1} = [\sigma_B ]_{AF}$. 
\endroster
\endproclaim
Let $\Aut(\FA,\DA)$
be the group of all automorphisms $\alpha$ of 
$\FA$ such that $\alpha(\DA) = \DA$.  
Denote by $\Inn(\FA,\DA)$ the subgroup of $\Aut(\FA,\DA)$ 
of inner automorphisms
on $\FA$.
We set 
$
\Out(\FA,\DA) 
$
the quotient group
$\Aut(\FA,\DA)/\Inn(\FA,\DA).$
By the same argument as Section 6, we have
\proclaim{Theorem 7.5}
There exist short exact sequeneces:
\roster
\item
$
1 
\longrightarrow Z_{\sigma_A}^1({\Cal U}(\DA))
\overset{\lambda}\to{\longrightarrow} \Aut(\FA,\DA)
\overset{\phi}\to{\longrightarrow} N([\sigma_A ]_{AF})
\longrightarrow 1
$
\item
$
1 
\longrightarrow B_{\sigma_A}^1({\Cal U}(\DA))
\overset{\lambda}\to{\longrightarrow} \Inn(\FA,\DA)
\overset{\phi}\to{\longrightarrow} [\sigma_A ]_{AF}
\longrightarrow 1,
$
\item
$
1 
\longrightarrow H_{\sigma_A}^1({\Cal U}(\DA))
\overset{\lambda}\to{\longrightarrow} \Out(\FA,\DA)
\overset{\phi}\to{\longrightarrow} N([\sigma_A ]_{AF}) / [\sigma_A ]_{AF}
\longrightarrow 1.
$
\endroster
They all split. 
Hence  
$
\Out(\FA,\DA)
$
is a  semi-direct product
$$
\Out(\FA,\DA)
= N([\sigma_A ]_{AF}) / [\sigma_A ]_{AF} \cdot H_{\sigma_A}^1({\Cal U}(\DA)).
$$ 
\endproclaim
where
$ N([\sigma_A ]_{AF})$ is the normalizer subgroup
of
$[\sigma_A ]_{AF}$
in 
$[[\sigma_A]]$.


\bigskip

{\it e-mail}: kengo{\@}yokohama-cu.ac.jp

\bye